\documentclass[a4paper,10pt]{article}
\usepackage{amsfonts}
\usepackage{amssymb}
\usepackage{graphicx}
\usepackage{mathtools}
\usepackage{skak}
\usepackage{MnSymbol}

\usepackage{amscd}
\usepackage{MnSymbol}

\usepackage{wasysym}

\usepackage{ifsym}

\makeatletter

\newcommand{\Rmnum}[1]{\expandafter\@slowromancap\romannumeral #1@}
\makeatother

\title{Hilbert Volume in Metric Spaces, Part 1.}

\author{Misha Gromov}

\begin{document}

\maketitle \tableofcontents

\begin{abstract} We introduce a notion of {\it Hilbertian $n$-volume} in metric spaces  with Besicovitch-type inequalities built-in into the definitions. The present Part 1 of     the article  is, for the most part,  dedicated to reformulation of known results in our terms with proofs being reduced to (almost) pure tautologies. If there is any novelty in the paper, this is in forging certain terminology which, ultimately, may turn useful  in  Alexandrov kind of approach to  singular spaces with positive scalar curvature \cite {bil}.

\end{abstract}

\section {  Partitions of Unity  $\mu=d_\mu p$ in Banach Spaces, $\tilde L_q$-Dilation $||min.\widetilde {dil}^\ast ||_{L_q}$
and  Hilbertian $\widetilde{Hilb}_{/id}$.}
We recall in this section a few (more or less)  standard definitions   
 and introduce  notations used throughout the paper.  
   
  \vspace {1mm}

   {\it  The      dilation} of a map between metric spaces,   $f:X\to Y$, is
   a function in two variables
     $x_1\neq x_2 \in X$, that is 
  $$dil_f(x_1,x_2)=dist_Y(f(x_1),f(x_2)/dist_X(x_1,x_2), \mbox {   $x_1\neq x_2 \in X$}.$$

{\it The   Lipschitz constant of $f$} is
  $$Lip(f)= \sup_{x_1\neq x_2 \in X}dil_f(x_1,x_2).$$

Equivalently,
 $$Lip(f)=\sup_{B\subset X} diam_Y(f(B))/diam_X(B),$$
where the supremum is taken over all bounded subsets in $X$. 
 
  \vspace {1mm}
  
 Let ${\cal L} = {\cal L} (Y)$ be the space of (usually Lipschitz) functions $l: Y\to \mathbb R$,
let ${\cal L}/const$ be the space of functions modulo additive
 constants, and let
$\mu=d_\mu l$ be a Borel measure on  ${\cal L}/const$.
Define the 
{\it $L_q$-dilation} of a map $f:X\to Y$, $1\leq q\leq \infty$,  as 
   $$||dil^\ast f||_{L_q(\mu)}=\left (\int _{{\cal L}/const}Lip^  q(l\circ f)d_\mu  l\right)^{1/q}.$$
  
\vspace {1mm}

{\it Example.}   Let $Y$ be the Euclidean space $\mathbb R^n$ and let $\mu$ be supported on the $n$ orthogonal projections (modulo constants)  of $\mathbb R^n$ onto
the coordinate axes with equal  weights 1 assigned to all projections. Then every isometric map $f: \mathbb R^m  \to\mathbb R^n$ satisfies
$$||dil^\ast f||_{L_2(\mu)} =\sqrt{m}.$$
 
 \vspace {1mm}
 
 In general, {\it an  axis} or a {\it (parametrized) straight line}  in a metric space $Y$, is an  isometric embedding  of the real line $\mathbb R$ into $X$, say $a  :\mathbb R\to Y$.  The image of   an axis  $\mathbb R=\mathbb R_a \subset Y$
 always admits a {\it $1$-Lipschitz  projector} of $Y$ onto it,  i.e. a map 
 $p: Y\to 
 \mathbb R_a \subset Y$ with $Lip(p)=1$ and such that $p^2=p$, i.e. 
  $p_{| \mathbb R_a} :  \mathbb R_a\to \mathbb R_a=id_{ \mathbb R_a}$.
  
  Such a projector is non-unique but in certain cases there are preferred ones.
  For example if $Y$ is a Banach space and $ \mathbb R_a\subset Y$ is a $1$-dimensional linear subspace, then it admits a  
  a {\it linear} projection $p: Y \to  \mathbb R_a$ by the {\it Hahn-Banach} theorem, where this $p$ is unique if the  Banach norm $||...||_Y$ is smooth.
  
There is a similar projector on an axis   $\mathbb  R_a \subset Y$ in an arbitrary metric space $Y$ called  {\it the Busemann function} $b_a :Y\to \mathbb  R= \mathbb  R_a \subset Y$ defined with $r\in  \mathbb R_a $ tending to $+\infty$ by 
  $$b_a(y)= \lim_{r\to  +\infty}\big (dist(a(r),y)-r)\big).$$ 
  
For instance, if $Y$ is a Banach space and the norm is smooth at the points $y\in    \mathbb R_a$ then $b_a$  coincides with the  Hahn-Banach projection.

  \vspace {1mm} 
{\it  An axial  projector}  in a Banach space $Y$ is a pair $(a,p)$, where
 $a: \mathbb R\to Y$ is a linear isometric imbedding $\mathbb R=\mathbb R_a \subset  Y$, and $p : Y\to \mathbb R_a=im_a=im_p\subset Y$ is  a linear  projector of $Y$ onto this axis with $Lip(p)=1$. 
 
 Observe that  $a$ is uniquely determined by $p$ up to $\pm$-sign; it will be sometimes denoted $a_p$ or  suppressed from the notation.

 Denote by ${\cal P}={\cal P}(Y)$  the space of axial  projectors $(a,p)$ in $Y$ and,
 given a  Borel measure  $\mu=d_\mu p$ on $\cal P$, define {\it  the $L_q$-dilation} of a map $f:X\to Y$, where $X$ is a metric space, by 
 
  $$  ||dil^\ast f||_{L_q(\mu)}=\left (\int _{\cal P}Lip^  q(p\circ f)d_\mu  p\right)^{1/q}, \mbox {  }  1\leq q \leq \infty.$$

 {\it  Identity Example.} The  $L_q$-dilation of the identity map $id_Y: Y\to Y$ is  obviously expressed by the total mass of $\mu$, 
 $$  ||dil^\ast f||_{L_q(\mu)}=\mu({\cal P})^{1/q}.$$
  
  \vspace {1mm}
  
{\it Maps  $\tau$,  ${\cal I}_\mu$ and the Norm $||\widetilde {dil}^\ast f||_{L_q}$.}   
 The space $Y$ tautologically  embeds into the space $\Phi({\cal P})$ of functions  ${\cal P}\to \mathbb R$ for $\tau: y\to \phi(p)=p(y)\in \mathbb R_p=\mathbb R$.  A
  Borel measure  $\mu$ on $\cal P$ gives the   $L_q$-norm to  the space    $ \Phi({\cal P})$ that is then denoted  $L_q({\cal P},\mu)$, where the $L_p$-norm of $\tau$ is bounded by $\mu({\cal P})^{1/q}$.

A measure $\mu$ on $\cal P$ also delivers a {\it  linear $1$-Lipschitz}
map ${\cal I}_\mu:\L_1({\cal P},\mu) \to Y$ that sends functions $\phi:(a,p) \mapsto \mathbb R=\mathbb R_p\subset Y$ to   their $\mu$-integrals in $Y$,
$${\cal I}_\mu: \phi(p)\mapsto \int_{\cal P}a(\phi(p)) d_\mu p.$$

Define {\it  $\tilde L_q$-dilation}
 $$||\widetilde {dil}^\ast f||_{L_q(\mu)} \leq  || {dil}^\ast f||_{L_q(\mu)}$$
via Lipschitz maps $\tilde f :X\to L_q({\cal P},\mu)$, that are $\cal P$-families of $\mathbb R$-valued functions  $\tilde f_p: X \to \mathbb R=\mathbb R _a\subset Y$, $(a,p)\in \cal P$,
  as  
$$||\widetilde {dil}^\ast f||_{L_q(\mu)}=\inf_{\tilde f} \left (\int_{\cal P}Lip^q(\tilde f_p) d_\mu p\right)^{1/q} $$
where the infimum is taken over all $\tilde f$  such that ${\cal I}_\mu\circ \tilde f=f$.

\vspace {1mm}

{\it An axial   partitions of unity in $Y$} is a measure $\mu$ on $\cal P$, such that
$$ {\cal I}_\mu\circ \tau=\int_{\cal P}p d_\mu p=id_Y\mbox { for the identity operator } id_Y: Y\to Y.$$

\vspace {1mm}

{\it Minimal $L_q$-Dilations.} Let
$$||min.dil^\ast f||_{L_q}=\inf_{\mu} ||dil^\ast f||_{L_q(\mu)}\mbox { and } ||min.\widetilde{dil}^\ast f||_{L_q}=\inf_{\mu} ||\widetilde{dil}^\ast f||_{L_q(\mu)}$$
where the infima are taken over all partitions of unity $\mu=dp_\mu$ in $Y$.

Notice that if $dim(Y)<\infty$ then both  minima are assumed by some partition of unity, say
$\mu_{min}(f)$ and $\mu_{\widetilde{min}}(f)$, such that
$$||min.dil^\ast f||_{L_q}= ||dil^\ast f||_{L_q(\mu_{min})}\mbox { and }
||min.\widetilde{dil}^\ast f||_{L_q}= ||\widetilde{dil}^\ast f||_{L_q(\mu_{\widetilde{min}})}.$$

\vspace {1mm}

{\it The Maps $\tau$ and ${\cal I}_\mu$ in Hilbert Spaces.}  If $Y$  is a {\it Hilbert  space}, then the tautological map $\tau:Y\to \Phi{(\cal P)}$ is {\it linear} and  if $\mu$ is a {\it partition of unity}, then $\tau:Y\to L_2({\cal P},\mu)$ is {\it isometric}  by the Pythagorean theorem.
 
The map ${\cal I}_\mu:\L_1({\cal P},\mu) \to Y$ extends to  ${\cal I}_\mu:\L_2({\cal P},\mu) \to Y$ that is,  in the Hilbertian case,  equals the adjoint to $\tau$, since the operators $a:\mathbb R\to Y$ and $p: Y\to \mathbb R_a=\mathbb R$ are mutually adjoint. Explicitly, 
$$\langle {\cal I}_\mu(\phi(p)), y\rangle_Y= \left\langle \int_{\cal P} a(\phi(p)) d_\mu p,y\right \rangle_Y=$$
$$= \int_{\cal P} a(\phi(p))\cdot p(y) d_\mu p=  \int_{\cal P} \phi(p)\cdot \tau(y) d_\mu p =\langle \phi(p), \tau(y)\rangle_{L_2({\cal P},\mu)}.$$

In other words,  ${\cal I}_\mu$ equals the orthogonal projection of ${L_2({\cal P},\mu)}$ onto $Y=\tau(Y) \subset {L_2({\cal P},\mu)}$. Consequently,  
{\it $$Lip({\cal I}_\mu)=1\mbox { for all partitions of unity $\mu$ in Hilbert spaces}.$$}
It follows, that every Lipschitz map $\tilde f$ from a metric space $X$ to  $L_2({\cal P},\mu)$, that is  a $\cal P$-family of $\mathbb R$-valued functions  $\tilde f_p: X \to \mathbb R=\mathbb R _a\subset Y$, $(a,p)\in \cal P$,  satisfies
$$ ||dil^\ast ({\cal I}_\mu\circ\tilde f)||_{L_2(\mu)}\leq\left  (\int Lip^2(\tilde f_p) d_\mu p\right )^{1/2}.$$

  {\it Partitions of Hilbertian Forms $h$ into Squares $l^2$.}   Let  $h$ be the Hilbertian quadratic form in $Y$.  Then the partition of unity condition on $\mu$ can be equivalently (and obviously) expressed in term of the integral
of the squares of linear functions (forms):
$$ \int _{\cal P}p d \mu_p =id_Y   \Leftrightarrow \int_{\cal L} l^2d \mu_l  =h,$$
where $l=l_p: Y\to \mathbb R=\mathbb R_p$ are linear functions corresponing to the projectors $p: Y \to \mathbb R_p$ and $d \mu_l$ is  the pushforward  of $\mu$  from $\cal P$ to to $\cal L$ under the map $p\mapsto l_p$.

If $X$ and $Y$ are Hilbert spaces and $f:X\to Y$ is a linear map, then $||dil^\ast f||_{L_2(\mu)}$ {\it does not depend} on the partition of unity $\mu$ in $Y$: it equals the trace  of the induced  quadratic form on $X$
 $$||dil^\ast f||_{L_2(\mu)}=  \int_{\cal L}trace_Y(l\circ f)^2 d \mu_l=trace_X  \int_{\cal L} (l\circ f)^2d \mu_l = trace_X f^\ast(h),$$
 and it is called   
    the {\it Hilbert -Schmidt} norm
 $||f||_{L_2}=||dil^\ast f||_{L_2(\mu)}$ of $f$.

Furthermore, since the map  ${\cal I}_\mu:L_2({\cal P},\mu)\to Y$ is $1$-Lipschitz, 
$$ ||\widetilde {dil}^\ast f||_{L_2(\mu)} =||dil^\ast f||_{L_2(\mu)}.$$

For instance,  isometric linear maps   $f:X\to Y$
 have 
 $$ ||\widetilde {dil}^\ast f||_{L_2(\mu)}= ||dil^\ast f||_{L_2(\mu)}=||f||_{L_2}=\sqrt{dim(X)}$$ 
for all partitions of unity $\mu$; in particular
$$ ||min.\widetilde {dil}^\ast id_Y||_{L_2}= ||min.dil^\ast id_Y||_{L_2}  =  ||id_Y||_{L_2}=\sqrt{dim(Y)}.$$ 
In general, Lipschitz maps $f$  from metric spaces $X$ to Hilbertian $Y$ satisfy 
$$ Lip(f)\leq   ||min.\widetilde {dil}^\ast f||_{L_2} \leq ||min.dil^\ast f||_{L_2}\leq  Lip(f) \cdot \sqrt{dim(Y)}.$$

This suggests the notation 
 $$ Hilb(f)=Hilb_{/id}(f)=Hilb_{/id_Y}(f) = ||min.dil^\ast f||_{L_2}/||min.dil^\ast id_Y||_{L_2}$$ 
and 
$$ \widetilde {Hilb}(f)= \widetilde {Hilb}_{/id}(f)= \widetilde {Hilb}_{/id_Y}(f) = ||min. \widetilde {dil}^\ast f||_{L_2}/||min. \widetilde {dil}^\ast id_Y||_{L_2}.$$ 
These definitions make sense for all
 arbitrary metric spaces $X$ and   Banach spaces $Y$, where, if $Y$ is Hilbertian, 
$$ \widetilde {Hilb}(f)\leq Hilb(f) \leq Lip (f) \leq  \sqrt{dim(Y)}\cdot   \widetilde{Hilb}(f). $$

\vspace {1mm}
Up to some  point,
we formulate and prove obvious  general properties of the "norms"
$||min.dil^\ast f||_{L_q}$ and  $||min.\widetilde {dil}^\ast f||_{L_q}$ for all Banach spaces $Y$ and all  $q$ -- this does not cost us anything, but
only maps into Hilbert spaces and   $q=2$ will be essential   for the present day geometric applications.

\vspace {1mm}

{\it Local Dilations $||dil^\ast f(x)||$, $Hilb(f(x))$ etc. in Finsler Spaces.} All of the above  dilation invariants  depends only on the metric but not on the linear structure in $Y$. This allows the following local versions of these dilations for maps into  Finsler, e.g Riemannian, manifolds $Y$.

Namely, let $Y=(Y, dist)$ be a metric space, where 
  each point $y\in Y$ admits a local Banach metric  in a neighbourhood $U_y \subset Y$ of $y$, say $dist_y$ for all $y\in Y$,
 such that 
 $$dist_y(y_1,y_2)/dist(y_1,y_2)\to 1 \mbox {  for } y_1,y_2\to y.$$
Define

 \hspace {2mm} $||{dil}^\ast f(x)||_{L_q(\mu)}$,  $||min. {dil}^\ast f(x)||_{L_q}$,   $|| \widetilde{dil}^\ast f(x)||_{L_q(\mu)}$, Hilb(f(x)),  etc.

\vspace {1mm}
 
 \hspace{-6mm}  for maps $f:X\to Y$
by restricting $f$ to small neighbourhoods $U_x\subset X$ for all $x\in X$, such
that $f(U_x)\subset U_{f(x)}\subset Y$, evaluating the corresponding dilations with respect to the local Banach metrics $dist_y$ and then taking the infimum over all neighbourhoods $U_x$ of a point  $ x\in X$ for all $x\in X$.

Equivalently,  in the case $Y$ is a Riemannian manifold, one can, by Nash' theorem,  isometrically immerse $Y\subset \mathbb R^N$ and define the local dilations
via  axial partitions of unity in  $\mathbb R^N$.

Notice, that unlike the Lipschitz constant, where $Lip(f)=\sup_{x\in X}Lip(f(x))$
for  maps from {\it length} spaces $X$, there is no(?) apparent passage from local to global $L_q$-dilations. 

\section {  Hadamard's  $Jac^{[n]} \leq \widetilde {Hilb}^n $
 and Inverse  Lipschitz Maps.}

   {\it  Hadamard's Inequality for $Jac^{[n]}$.} Given a map between Hilbert spaces, say $D:A\to B$, denote by  $Jac^{[n]}(D)$ the norm of the $n$-th exterior power of $D$, that is the supremum of the {\it absolute values} of the  Jacobians of $D$ on all $n$-dimensional subspaces $A'\subset A$,
  $$Jac^{[n]}(D)=\sup_{dim(A')=n} Jac (D|A')\mbox { for } D|A': A'\to im_{A'} \subset B.$$

 If $f:X\to Y$  is a locally Lipschitz, e.g. $C^1$-smooth, map between Riemannian manifolds 
 then the differential $Df(x):T_x(X) \to T_{f(x)}(Y)$ exists for almost all $x\in X$
by {\it Rademacher-Stepanov theorem} and the Jacobian of $f$ is defined as   
 $$Jac^{[n]}(f(x))=Jac^{[n]}(Df(x)) \mbox  {  and } Jac^{[n]}(f) =\sup_ {x} Jac^{[n]}(Df(x)),$$
 where $\sup_x$  refers to {\it almost all} $x\in X$.
 
 \vspace{1mm}

{\it   Hadamard's Inequality.}  Let $X$ and $Y$ be  Riemannian manifolds  and $f: X\to Y$ be a Lipschitz map.

{\it Then
 $$Jac^{[n]}(f(x)) \leq  ||n^{-1/2}min.\widetilde{dil}^\ast f(x)||^n_{L_2}\mbox {  for $n=dim(Y)$ and almost all $x\in X$}. $$
In particular, Lipschitz maps $f$ from $n$-manifolds to $\mathbb R^n$ satisfy
$$Jac^{[n]}(f) \leq \widetilde {Hilb}^n(f).$$}

{\it Proof.} In fact, if
$$f={\cal I}_\mu \circ  \tilde f =\int_{\cal P} a_p(\tilde f) d_\mu p \mbox { for $\tilde f: X\to L_2({\cal P},\mu)$},$$
 then
 $$ ||D f(x)||_{L_2}\leq ||D\tilde f(x)||_{L_2} \leq ||min.\widetilde{dil}^\ast f(x)||_{L_2}\mbox   { for almost  all } x \in X,$$
while 
 $$Jac^{[n]}(D) \leq n^{-1/2} ||D||^n_{L_2},$$
 since the discriminants of  quadratic forms (induced by $D$ from the Hilbert form in target Hilbert spaces) are bounded by their traces via the
 arithmetic/geometric mean inequality for the eigenvalues of these forms.

\vspace {1mm}

However trivial, the inequality   $Jac^{[n]}(f) \leq Hilb^n(f)$  is significantly  sharper than mere $ |Jac^{[n]}|(f)\leq Lip^n(f)$. 
 
 For instance, let a  smooth map $f$ from a Riemannian $n$-manifold $X$ to 
 $Y=\mathbb R^n$ be given by $n$ functions $f_1,...f_n$.
 Then  $Lip(f)$ may be as big as $(\sum_i Lip^2(f_i))^{1/2}$, while the Hadamard's inequality, applied  to the Hilbert form
$\sum_{i=1}^n  \lambda_ib_i^2 $ for $\lambda_i=(Lip(f_i))^{-1/2}$ on $Y$,
says that 
 $$ |Jac^{[n]}(f)|\leq \prod _{i=1,...,n}Lip(f_i).$$
(This is, of course,  obvious anyway.) 

\vspace {1mm}

 {\it Sharpness of Hadamard.}  Hadamard's   inequality  bounds the volume of the image of a  Lipschitz map between $n$-dimensional Riemannian manifolds, $f:X\to Y$, 
    by $$ vol_n(f(X))\leq \int _XHilb^n(Df(x))dx\leq
  vol_n(X)\cdot  \sup_xHilb^n(Df(x)),$$
since
   $$vol_n(f(X))\leq \int _XJac^{[n]}(f(x))dx \leq  vol_n(X)\cdot   \sup_xJac^{[n]}(Df(x)),$$
 where  the two suprema are taken over {\it almost all} $x\in X$.

  \vspace {1mm}
  Since the algebraic inequality  $Jac^{[n]}(D) \leq Hilb^n(D)$
  is sharp for linear maps $D:\mathbb R^n\to \mathbb R^n$, i.e. the equality   $Jac^{[n]}(D) = Hilb^n(D)$  implies  that $D$ is a homothety, 
  the corresponding integral inequalities are also sharp:     
   \vspace {1mm}
 
 \textbf {A.} {\it  If  $$vol_n(f(X))= \int _XHilb^n(Df(x))dx<\infty\mbox   { and  }dim (X)=dim(Y)\geq 2,$$
 then the map $f$ is  a conformal diffeomorphism on its image. } 

(This is not  used in the sequel.)

  \vspace {1mm}

 \textbf {B.} {\it If 
  $$   \sup_xHilb^n(Df(x))=1 \mbox { and } vol_n(f(X))= vol_n(X)<\infty,$$ 
  then $f$ is an isometric diffeomorphism onto its image.}
 
 \vspace {1mm}
 
{\it Proof.}  The  equality $vol_n(f(X))= \int _X Jac^{[n]}(Df(x))dx$ implies that the local topological  degree of   $f$,   either equals  $+1$    at almost  all $x\in X$   or it is  almost everywhere $-1$. If such an $f$ is conformal {\it almost everywhere} and $n\geq 2$, then it is conformal {\it everywhere}. For example, if $n=2$
 this is seen with the  Cauchy integral formula, while  the case $n\geq 3$ is (essentially) trivial because of the Liouville theorem.
 
  The remaining case of  \textbf {B} for
 $n=1$ is, of course, obvious but it points out to  the following extension  to
more general Lipschitz maps.

 \vspace {1mm}

 \textbf {C.}  {\it If a Lipschitz map $f:X\to Y$ between $n$-dimensional Riemannian manifolds without boundaries  satisfies  
 $$vol_n(f(X))= \int _XJac^{[n]}(f(x))dx <\infty$$
 and  
  $$ vol_n(f^{-1}(U))\leq C\cdot vol(U)\mbox { for all  open  $U\subset Y$ and a constant $C<\infty$.}$$
 Then $f$ is a homeomorphism onto its image $im_f=f(X)\subset Y$ and the inverse map
 $f^{-1} : im_f\to X$ satisfies
 $$\sup_{y\in im_f}Lip_y(f^{-1})\leq C\cdot Jac^{[n-1]}(f) \leq C\cdot Lip^{n-1}(f).$$}

This is well known with   the proof, probably,  buried somewhere in \cite {fed}.  
  Here is how it goes.

Since  $ vol_n(f^{-1}(U))\leq C\cdot vol(U)$, the pullbacks  $f^{-1}(y)\subset X$ are  {\it zero dimensional} subsets in $X$ for all  $y\in im_f$;  in fact, this follows from   the inequality
$vol_n (f^{-1}(B_y(\varepsilon)))=o(\varepsilon^{n-1})$ for small $\varepsilon$-balls in $Y$.

Since  $dim(f^{-1}(y)) \leq 0$,  every point $x\in f^{-1}(y) $ admits an  arbitrarily small neighbourhood, say 
  $B_x\subset X$, such that 
 the  the boundary   $\partial B_x \subset X$  does not intersect $f^{-1}(y)$; hence, the local  topological degree of $f$ is defined at $x$. 
 
 This degree must be {\it non-zero}, otherwise, all points $x'$ close to $x$ would have "partners" $x''$ such that $f(x')=f(x)$ that is incompatible with   $vol_n(f(X))= \int _XJac^{[n]}(f(x))dx$. 
 
  Therefore, the image of every  neighbourhood of $x\in X$ contains a neighbourhood of $y=f(x)\in Y$ for all $x\in X$. Thus, $f$ is what is called an {\it open map}. 
  
  Besides, again because of  $vol_n(f(X))= \int _XJac^{[n]}(f(x))dx$, the map  
$f$ is one-to-one on a dense subset (of full measure) in $X$, and, obviously,
   
  \hspace {16mm} densely one to one + open $ \Rightarrow$ one-to-one.

  Granted one-to-one, take an  $\varepsilon$-narrow cylinder $V_\varepsilon \subset im_f$ around a distance minimizing  geodesic segment $[y_1,y_2] \subset im_f$, denote  by $\partial_1 V_\varepsilon\ni y_1$ and    $\partial_2 V_\varepsilon\ni y_2$. the top and the bottom of this cylinder and  let $\partial_{lat}V_\varepsilon$ be the remaining (lateral) part of its bounadry,  
  $$\partial_{lat}V_\varepsilon =\partial V_\varepsilon \setminus( \partial_1V_\varepsilon \bigcup \partial_2V_\varepsilon).$$ 
Observe that the relative homology group $H_{n-1}(V_\varepsilon, \partial_{lat}V_\varepsilon)$ is free cyclic and let 
$\kappa_{n-1} \varepsilon^{n-1}$ denote the minimum of the $(n-1)$-volumes of relative $(n-1)$-cycles representing {\it non-zero} classes in  $H_{n-1}(V_\varepsilon, \partial_{lat}V_\varepsilon)$.

Observe that $\kappa_{n-1}$ is close to the volume of the Euclidean $(n-1)$-ball 
and 
$$vol_n(V_\varepsilon)=\kappa_{n-1} \varepsilon^{n-1}\cdot dist (y_1,y_2) + o (\varepsilon^{n-1}).$$

Take the pull-backs $\tilde V_\varepsilon =f^{-1}(V_\varepsilon)
\subset X$ and $\partial_{lat}\tilde V_\varepsilon=f^{-1}(\partial_{lat}V_\varepsilon) $ and let $\tilde \kappa_{n-1} \varepsilon^{n-1}$
 be the minimum of the $(n-1)$-volumes of relative $(n-1)$-cycles representing {\it non-zero} classes in  $H_{n-1}(\tilde V_\varepsilon, \partial_{lat}\tilde V_\varepsilon)$.

Clearly,
$$ \kappa_{n-1} \leq   Jac^{[n-1]}(f)\cdot \tilde \kappa_{n-1}.$$

On the other hand,
 $$vol_n(\tilde V_\varepsilon))\geq   dist _X( \partial_1\tilde V_\varepsilon,
 \partial_2\tilde V_\varepsilon)\cdot \tilde \kappa_{n-1} \varepsilon^{n-1}$$
for   $\partial_1\tilde V_\varepsilon, \partial_2\tilde V_\varepsilon \subset X$ being the pullbacks of the top and the bottom of the cylinder  $V_\varepsilon \subset Y$. 

The proof of this is standard:  map $\tilde V_\varepsilon$ to $\mathbb R_+$ by 
 $d:\tilde v\mapsto dist_X(\tilde v,  \partial_1\tilde V_\varepsilon)$,
observe that the pullbacks $d^{-1}(r) \subset \tilde V_\varepsilon$ support non non-zero  classes
in $H_n(\tilde V_\varepsilon,  \partial_{lat}\tilde V_\varepsilon)$ for $0\leq r\leq
dist_X (  \partial_1\tilde V_\varepsilon,  \partial_2\tilde V_\varepsilon)$ and apply the classical (and obvious) {\it  coarea inequality} (that happens to be equality in the present case and that is extended in section 8 to more general setting of Hilbert volumes),
$$vol_n(\tilde V_\varepsilon)\geq \int_ {\mathbb R_+} vol_{n-1}(d^{-1}(r)) dr.$$
Then the proof follows with $\varepsilon \to 0$.

 \vspace {1mm}

\textbf {D.} {\it Remark.} The statement of \textbf {C} can be obviously reformulated without any reference to differentiability of $f$ and  
     the above argument still applies as it  does not use the  Rademacher-Stepanov theorem.
     
Besides, one does not truly need $f$ to be Lipschitz:  what is essential is
 a bound on the  $(n-1)$-Jacobian $Jac^{[n-1]}(f)$ that can be   defined as 
 $$Jac^{[n-1]}(f)= \sup_H  vol_{n-1}(f(H))/ vol_{n-1}(H)$$   
 for  all rectifiable hypersurfaces $H\subset X$.
 
 Moreover, $X$ and $Y$ do not have to  carry  any metrics, what is needed are measures on $X$ and $Y$ and the notion of "$(n-1)$-volume" for hypersurfaces.   
  In fact, the distance can be {\it derived}  from these data via 
 the {\it Almgren-Besicovitch-Derrick inequality} (see section 9).     
    
     \vspace {1mm}
 
 { \it Warnings.} (a)   A Lipschitz map $f:X\to Y$ where $Df(x):T_x(X)\to T_{f(x)}(Y)$ is an isometry almost everywhere does not have to be one-to-one.
 In fact,  paradoxically, every Riemannian $n$-manifold  admits a $1$-Lipschitz map into
 $\mathbb R^n$, that preserves the lengths of {\it all} rectifiable curves in $X$,   see  2.4.11 in  \cite{pdr}.
  
   (b) There are lots of   $C^\infty$-smooth maps $f: X\to Y$ such that $vol_n(f(X))= \int _XJac^{[n]}(f(x))dx$,  where $Jac^{[n]}(f)\neq 0$ almost everywhere; yet, where the  unions of those pullbacks $f^{-1}(y) \subset X$ that are diffeomorphic to the  $(n-1)$-ball, are {\it dense} in $X$.

  \section { Controlled $\tilde L_q$-Dilation Extensions,  Banach Straightening and Hilbertian Volume Domination.}

If $ X\supset  X_0  \overset {f_0}\to Y$ is a  partially defined  $\lambda$-Lipschitz map between metric spaces, then it admits an extension to a  $\lambda$-Lipschitz
map $f: X\to Y$  only for rather exceptional spaces $Y$. For instance, this is, obviously, possible if   $Y$ is a tree,
e.g. $Y=\mathbb R$,   a (finite or infinite measurable)  Cartesian product of trees with the  $sup$-metric,  e.g. the $L_\infty$-space over a measure space, such as    $ L_\infty({\cal P}, \mu)$,
or else,
 an {\it Uryson's universal space}.

But such  extensions do not exist, in general,  not even for $Y=\mathbb R^n$ if $n\geq 2$,
except for particular domains $X$, such as {\it  Alexandrov's spaces} with curvature $\geq 0$  \cite {lan-sch}, while for general $X$, only  an extension $f:X\to Y=\mathbb  R^n$ with $Lip(f)\leq n^{1/2}Lip(f_0)$ is (obviously) possible.

However, Lipschitz extensions of maps $\tilde f_0: X_0\to   L_\infty({\cal P}, \mu)$ to $\tilde f :  X\to   L_\infty({\cal P}, \mu)$  trivially  yield a  control over $\tilde L_q$-dilation of  maps $f={\cal I}_\mu\circ \tilde f :X\to Y$ as follows.

\vspace {1mm} 

Let $X$ be  a metric space and  $Y$ be a  Banach space with  an axial partition of unity $\mu$.
 
\vspace {1mm}

{\it  Then every partial map,   $X\supset X_0\overset {f_0} \to Y$, extends to a map $f:X\to Y$,    such that
  $$||\widetilde{dil}^\ast f||_{L_q(  \mu)} \leq  ||\widetilde{dil}^\ast f_0||_{L_q(  \mu)}\leq  ||{dil}^\ast f_0||_{L_q(  \mu)}\leq$$
  $$ \leq ||dil^\ast f_0||_{L_\infty( \mu )} \cdot \big (\mu({\cal P})\big )^{1/q}\leq Lip(f_0)\cdot \big (\mu({\cal P})\big )^{1/q}$$   
  \hspace {50mm} for all $   1\leq q\leq \infty.$}

\vspace {1mm}

{\it Remark.} The only point here which needs a minor attentions (and which is unneeded for  our applications) is to make sure that the same
 extension $\tilde f$  of  $\tilde f_0$  serves the spaces $L_q({\cal P},\mu)$ simultaneously for all $q$.
 
  This is achieved by  making Lipschitz extensions 
  $\tilde f_p: X \to \mathbb R=\mathbb R_p\subset Y $, $p\in \cal P$, of the corresponding functions    $(\tilde f_0)_p: X_0 \to \mathbb R=\mathbb R_p\subset Y $ depend {\it measurably} on $p$, e.g. by taking, for each $p\in \cal P$, a {\it minimal} such $\tilde f_p$. Namely, in general,
  given a partial defined function   $X\supset X_0\overset {\phi_0} \to \mathbb R$,
  its 
 minimal Lipschitz extension is
$$\phi(x)=\sup_{x_0\in X_0} (\phi_0(x_0)-dist_X(x,x_0)).$$   
For example, if $X_0$ is isometric to an axis $\mathbb R\subset X$ and $X_0\to
\mathbb R=X_0$ is the identity map, this minimal extension equals the Busemann function $X\to \mathbb R$.

 \vspace {1mm}

{\it Banach Straightening  at Infinity.}  Define the  {\it Lipschitz constant at infinity}  of a map $f:X\to Y$ between metric spaces by
 $$Lip_\infty(f)=\limsup_{dist_X(x_1,x_2)\to \infty}dil_f(x_1,x_2). $$

Equivalently, this can be defined with restrictions of $f$ to 
 {\it $D$-separated  nets} $X_D \subset X$, that are  subsets in $X$ with distances between all pairs of points $\geq D$, as follows,
$$Lip_\infty ( f)= \limsup_{X_{ D\to \infty}}Lip(f_{|X_D}).$$

Similarly, define the $\tilde L_q$-dilation at infinity, denoted 
$||\widetilde {dil}_\infty^\ast f||_{L_q(\mu)}$, by restricting $f$ to     $D$-separated nets
and taking   $limsup$ of $ ||\widetilde {dil}^\ast f_ {|X_D}||_{L_q(\mu)}$   over $X_D$   with $D\to \infty$.

 \vspace {1mm}
  
Then the above extension of  maps  $f|X_D:X_D\to Y$  from  $X_D\subset X$ to all of $X$ implies the following

 \vspace {1mm}

  {\it Banach $\varepsilon$-Straightening.}   Let $f:X\to Y$ be a map from a metric space $X$ to a Banach space $Y$ with a partition of unity $\mu$ in $Y$.
  
  \vspace {1mm}

{\it Then, for every $\varepsilon>0$, there exists a map $f_\varepsilon:X\to Y$   within finite distance from $f$,
i.e.
$$||f_\varepsilon(x)-f(x)||_Y \leq r=r_\varepsilon <\infty\mbox  { for  all } x\in X,$$
and such that
  $$ ||\widetilde {dil}^\ast f||_{L_q(\mu)}\leq ( ||dil_\infty^\ast f||_{L_\infty(\mu)}+\varepsilon)\cdot \mu({\cal P})^{1/q}\mbox { for all }  q\in [1,\infty].$$}

{{\it On $\varepsilon \to 0$}. One can, in many cases, pass to the limit map $f_{\varepsilon\to 0}:X\to Y$ but this  $f_{\varepsilon\to 0}$ may be far away from $f$, e.g. it may be  constant
for an $f$  {\it isometric at infinity}, i.e. where  $ \lim_{D\to \infty} dil_{f|X_D}=1$.

\vspace {1mm}

{\it Hilbert-Hadamard Volume Domination.}
Let  $X$ be an $n$-dimensional Riemannian manifold, let $\mu$ be a partition of unity in the Euclidean/Hilbertian  space   $Y=\mathbb R^n$ and  
let  $f: X\to  Y$  a  proper  (pullbacks of compact sets are compact) continuous  map. 

Let $\Omega_R\subset \mathbb R^n$, $R\to \infty$, be bounded domains, such that 
$vol_n(\Omega_R)= const\cdot  R^n$ and such that for every $r_0$ the volumes
of the $r_0$-neighbourhoods of their boundaries satisfy
$$ vol_n(U_{r_0}(\partial \Omega_R))=o(R^n)\mbox  { for } R\to \infty,$$
e.g. $\Omega_R$ are Euclidean  $R$-balls. 

\vspace {1mm}
 {\it  If $f$ has non-zero topological degree, then  the Euclidean volumes  of  the domains
      $\Omega_R  \subset \mathbb R^n$ are "asymptotically  sharply bounded" by the Riemannian volumes of their $ f$-pullbacks as follows 
      $$  vol_n\big( f^{-1}(  \Omega_R )\big)\geq ||dil_\infty^\ast||^n_{L_\infty(\mu)} \cdot vol_n(  \Omega_R )  -o(R^n)\mbox { for } R\to \infty.$$}

\vspace {1mm}

{\it Proof.} Since $deg(f)\neq 0$ and $dist(f, f_\varepsilon)< r_\varepsilon$,
the images $f_\varepsilon (f^{-1} (  \Omega_R )) \subset \mathbb R^n$ contain  
   $ [ \Omega_R-r_\varepsilon]=_{def}\Omega_R\setminus U_{r_\varepsilon}(\partial \Omega_R)$;  hence,
$$vol_n(f^{-1} (  \Omega_R )) \geq vol_n(f^{-1}_\varepsilon  [ \Omega_R-r_\varepsilon]).$$

 On the other hand,  by the above, $$\widetilde{ Hilb}(f_\varepsilon)=n^{-1/2} ||\widetilde{dil}^\ast f_\varepsilon||_{L_2(\mu)}   \leq  ||dil_\infty^\ast f||_{L_\infty(\mu)} +\varepsilon$$
while the Jacobian of $f_\varepsilon$ is bounded by  Hadanmard's inequality 
$Jac^{[n]} (f_\varepsilon) \leq \widetilde {Hilb}^n(f_\varepsilon)$; hence, 
 $$ vol_n(f^{-1}_\varepsilon  [ \Omega_R-r_\varepsilon])\geq
   \widetilde {Hilb}^{-n}(f_\varepsilon)\cdot vol_n[  \Omega_R-r_\varepsilon]$$
$$\geq  ||dil_\infty^\ast f||^{-n}_{L_\infty(\mu)} \cdot vol_n(  \Omega_R ) -\varepsilon\cdot const\cdot R^n -o(R^n) \mbox { for all  $\varepsilon >0$}.$$
QED.

\vspace {1mm}

 {\it Local $\tilde L_q$-Extensions.} Let $Y$ be a smooth manifold with a continuous   Riemannian metric and  $X\supset X_0\overset {f_0} \to Y$ a partial Lipschitz map.
 Suppose that that 
 $$||min\widetilde {dil}^\ast f_0(x)||_{ L_q} < \phi(x)\mbox {  for  all $x\in X_0$} $$ where $\phi$ is a  continuous  function  on $X$.

\vspace {1mm}

{\it Then $f_0$ extends to a  neighbourhood $X_1\supset X_0$ by a map 
 $X\supset X_1\overset {f} \to Y$, such that 
 $$||min\widetilde {dil}^\ast f(x)||_{ L_q} < \phi(x)\mbox {  for  all $x\in X_0$}. $$}

{\it Proof.}  Combine the above extension with the ordinary Lipschitz partition of unity in $X$ by the following trivial argument.

First, let $Y=\mathbb R^N$ and  let $U_i\subset X$, $i\in I$, be a locally finite covering of $X$ such that 

 $$||min\widetilde {dil}^\ast f_{0|U_i \cap X_0}||_{ L_q} < \phi(x)\mbox  {  for  all $x\in U_i$}.$$
Extend the maps $ f_{0|U_i \cap X_0}: U_i \cap X_0\to \mathbb R^N$  to $f_i:U_i\to \mathbb R^N$, such that 
$$||min\widetilde {dil}^\ast f_i||_{ L_q} \leq ||min\widetilde {dil}^\ast f_{0|U_i \cap X_0}||_{ L_q}< \phi(x)\mbox  {  for  all $x\in U_i$}.$$
Let $\phi_i:U_i\to\mathbb R_+$ be Lipschitz functions with supports strictly inside
$U_i$ and  such that $\sum_{i\in I}$=1. Then the map
$$F=\sum_{i\in I}\phi_i\cdot f_i\mbox {
satisfies }
||min\widetilde {dil}^\ast F(x)||_{ L_q} < \phi(x)$$
in some neighbourhood $X_1\subset X$ of $X_0$ by an obvious computation.

Finally, if $Y$ is $C^0$-Riemannian, isometrically $C^1$-embed $Y\subset \mathbb R^N$
and observe that some neighbourhood $Y_1\subset \mathbb R^N$ of $Y$
admits a Lipschitz projection $\pi:Y_1\to Y$ such that $Lip(\pi(y))=1$ for all 
$y\in Y$. (If $Y\subset \mathbb R^N$ is $C^2$-smooth, one can use the normal projection to $Y$, that does not always  exist for $C^1$-submanifolds.)

Compose the above $F$ with this $\pi$ and thus, obtain the required extension $f=\pi\circ F:X\to Y$ of $f_0$ from $X_0$  to $X_1\supset X_0$. QED.

\vspace {1mm}
{\it Question on $L_q$-Dilation at Infinity.}  Is there a counterpart to the above with $||dil^\ast f||_{L_q}$ instead of $||  \widetilde {dil}^\ast f||_{L_q}$?
 
 Namely, what is the minimal (or, rather, infimal) constant $C=C(Y,q, q')$, such that  every $f:X\to Y$ admits an
$f':X\to Y$ within bounded distance from $f$, such that
$||dil^\ast f'||_{L_{q'}(\mu)}\leq C\cdot ||dil_\infty^\ast f||_{L_q(\mu)}$?

\section {John's $h_\diamond$ and   John's  Ellipsoid in Banach Spaces with  a Riemannian Corollary.}  

  So far, everything was  compiled of definitions  + trivial generalities. Now comes something deceptively simple but more substantial.  
 
 \vspace {1mm}
  
  \textbf { Fritz John's Ellipsoid Theorem.}   
 {\it An $n$-dimensional Banach space  $Y=(Y,||...||_Y)$ admits a unique Hilbert quadratic form $h_\diamond$, such that 
  $$||...||_{h_\diamond}\geq ||...||_Y,$$
 and such that   the identity map $id_\diamond=id: Y=(Y,||...||_Y)\to Y=Y_\diamond=(Y,||...||_{h_\diamond})$ satisfies
$$ Hilb_{h_\diamond} =Hilb_{/Y_\diamond}(id_\diamond)\leq 1,$$
for  $Hilb_{/Y_\diamond}(id_\diamond)=n^{-1/2}||min. dil^\ast id_\diamond||_{L_2}$.}

\vspace {1mm}

{\it  Proof.} 
Let $H_+(1)$ be the subset in the space of  positive semidefinite (Hilbertian) quadratic forms  $h$  on $Y$ such that  $Hilb_{h}\leq 1$, i.e such that the identity map $id=id_h: (Y, ||...||_Y) \to (Y,h)$.
 satisfies  $n^{-1/2}||min. dil^\ast id_h||_{L_2} \leq 1$

 Observe that  this   $H_+(1)$ is a {\it convex} susbset in the linear space of all quadratic forms on $Y$ by the definition of 
 $||min. dil^\ast||_{L_2}$ via partitions of Hilbertian forms $h$ into  squares $l^2$ (see section 1).

\vspace {1mm}
 
$(\diamond)$  Let $h_\diamond\in H_+(1)$  {\it maximize} the  Hilbertian (Euclidean) Haar measure of $Y $ associated  to it among all  $h \in  H_+(1) $,  i.e.
the Jacobian  of the identity map $id : (Y,h_\diamond) \to (Y,h)$, denoted  $Jac^{[n]}_{ h_\diamond \to h}$, satisfies
$Jac^{[n]}_{ h_\diamond \to h} \leq 1$ for all $h\in H_+(1)$.
Then 
$$||...||_{h_\diamond} \geq ||...||_Y.$$

Indeed, 
 assume otherwise,  let $l$ be  a linear form  on $Y$ such that $||l||^2_Y =n=dim(Y) $ and $||l||^2_ {h_\diamond}=n/c$ for $c>1$ and 
let
$h_\varepsilon= (1-\varepsilon)h_\diamond +\varepsilon l^2$, $0\leq \varepsilon  \leq 1$.

Clearly, 
$$Jac^{[n]}_{ h_\diamond \to h_\varepsilon}=\left ( (1-\varepsilon)^{(n-1)}\cdot( (1-\varepsilon) +cn\varepsilon))\right)^{1/2},  $$
and
$$ \log Jac^{[n]}_{h_\diamond \to h_\varepsilon} =\frac{1}{2}\left ( -(n-1)\varepsilon +c(n-1)\varepsilon\right) -o(\varepsilon)>0\mbox {  for   }\varepsilon\to 0;$$
hence, the form $h_\varepsilon$ has greater Haar measure than $h_\diamond$ for small
$\varepsilon >0$.

Since $n\cdot l^2 \in H_+(1)$, the form $ h_\varepsilon$ also lies in  $ H_+(1)$ by convexity of  $H_+(1)$; this contradicts  the extremality of $h_\diamond$ and   uniqueness of $h_\diamond$ follows by  the same argument. QED.

\vspace {1mm}

{\it John's $\mu_\diamond$.} The partition of unity $\mu_\diamond=\mu_{min}$ in the Hilbert space  $(Y, ||...||_{h_\diamond})$ for which  $||dil^\ast id_\diamond||_{L_2}=\sqrt{n}$ is not, in general unique. However, the inequalities  $||...||_{h_\diamond}\geq ||...||_Y$ and 
$ Hilb_{h_\diamond}\leq 1$ trivially imply that

\vspace {1mm}

{\it every projector $p \in  {\cal P}(Y, ||...||_{h_\diamond})$  from the support of some  $\mu_\diamond$  lies in    ${\cal P}(Y, ||...||_Y)$, i.e. $||p||_Y= 1$, and the norm 
$||...||_{h_\diamond}$ equals $||...||_Y$ on the axes $\mathbb R_p=im_p\subset Y$ of all $p\in supp(\mu_\diamond)$.}

\vspace {1mm} 
 {\it Remark.} The above shows that the total masses  of John's partitions of unity
 $\mu_\diamond$  in $Y$ satisfy $\mu_\diamond({\cal P)}=n$=dim $Y$. In fact, the existence of a partition of unity in $Y$ with the total mass $\leq n$ is equivalent (by a trivial argument) to the full John's theorem.

\vspace {2 mm}

 {\it Riemannian  Lower Volume Bound.} Let $Y$ be an $n$-dimensional Banach space, let  $\Omega_R\subset  Y$, $R\to \infty$, be bounded domains, as in 
the   Hilbert-Hadamard volume domination  from section 3, i.e. 
 such that their $n$-volumes with respect to John's Hilbert/Euclidean metric $h_\diamond$ satisfy
$vol_\diamond(\Omega_R)= const\cdot  R^n$ and such that for every $r_0$ the volumes
of the $r_0$-neighbourhoods of their boundaries satisfy
$$ vol_\diamond(U_{r_0}(\partial \Omega_R))=o(R^n)\mbox  { for } R\to \infty,$$
e.g. $\Omega_R$ are the   $R$-balls $\{||y||\leq R\} \subset Y$.

\vspace {1mm}

  {\it   Asymptotic Banach-John Volume Inequality of  Burago-Ivanov.}  Let $X$ be an $n$-dimensional Riemannian manifold, let    $f:X\to Y$ be a proper (pullbacks of compact set are compact)
map of non-zero degree and let $Lip_\infty(f)$ be the Lipschitz constant
of $f$ at infinity with respect to the metric $||y_1-y_2||_Y$ in $Y$.
 
 \vspace {1mm}
 
  {\it Then the (John's Euclidean)  $h_\diamond$-volumes  of       $  \Omega_R  \subset  Y$ are "asymptotically  sharply bounded" by the Riemannian volumes of their $f$-pullbacks by  
$$  vol_n(f^{-1}(  \Omega_R ))\geq Lip^{-n}_\infty(f) \cdot vol_\diamond(  \Omega_R )  -o(R^n)\mbox { for } R\to \infty.$$
Consequently, since $||...||_{h_\diamond}\geq ||...||_Y$, the 
asymptotic volume growth of $R$-balls in $X$ around any given point $x_0\in X$ is asymptotically sharply bounded by the growth
of the Euclidean balls,
$$\liminf_{R\to \infty} vol(B_X (R))/ vol(B_{\mathbb R^n} (R))\geq 1.$$}

  \vspace {1mm}
 {\it Proof.}  Since John's  $\mu_\diamond$ serves as partition of unity in  $(Y, ||...||_Y)$ as well as in  $(Y,||...||_{h_\diamond})$, the map $ f :X\to  \mathbb R^n=(Y, ||...||_{h_\diamond})$  
 satisfies 
 $$||dil_{h_\diamond}^\ast f||_{L_\infty(\mu)} \leq  Lip_\infty(f)_{||...||_Y},$$
 where -- this is the main point --  the dilation $ ||dil^\ast f||$ is measured  with the
 (Euclidean)   $h_\diamond$-metric in $Y$, while   $Lip_\infty(f)$ is evaluated with the original Banach metric. Hence the  Hilbert-Hadamard volume domination applies and the proof follows.

\section {  Co-Lipschitz at Infinity,
Federer-Whitney Metric Descendants   and  Abelian Volume Growth Inequality of Burago-Ivanov.}

 {\it  Federer-Whitney Theorem.}  
  Define the {\it co-Lipschitz constant at infinity}   of a map between metric spaces, $f:X\to Y$,  by  
 $$coLip_\infty(f)=  \limsup_{diam_Y(B)\to \infty} diam_X(f^{-1}(B))/diam_Y(B),$$
   where $B\subset Y$ run over all {\it bounded} subsets in $Y$.
   
 Observe that if   
$$\sup_{y\in Y}  diam_X(f^{-1}(y)) <\infty,$$ 
(which is weaker than $coLip_\infty(f)<\infty$)
then  $$coLip_\infty(f)=\big (\liminf_{dist(x_1,x_2)\to \infty} dil_f(x_1,x_2)\big )^{-1}.$$

\vspace {1mm}
A map between metric spaces, $f:X\to Y$,  is called  an   {\it isometry at infinity}
 if  the image of $f$ intersect all $R$-balls in $Y$ for $R\geq R_0<\infty$ and
 $$Lip_\infty(f)=coLip_\infty(f)=1;$$
this is equivalent to    
   $$dil_f(x_1,x_2)=\frac{dist_X(x_1,x_2)}{dist_Y(f(x_1)f(x_2))}\to 1\mbox {  for } dist_X(x_1,x_2)\to \infty.  $$

 \vspace {1mm}
 
 Let  $X$ be a locally compact metric space,  $Y$  an $\mathbb R$-linear $n$-space, where both spaces are acted upon by a locally compact group $\Gamma$, such that the action of $\Gamma$ on $X$ is {\it isometric} and the action on $Y$ is {\it affine}. 
 
 Let this  affine action be {\it co-compact quasi-translational}, i.e.  $Y/\Gamma$ is compact  and $\Gamma$  admits a {\it co-compact}  subgroup   that acts on $Y$ by {\it parallel translations}.

\vspace {1mm} 

{\it  Federer-Whitney  $f$-Descendant Metric.}     Let $f : X \to Y$ be a proper   continuous $\Gamma$-equivariant map. 

\vspace {1mm} 

{\it Then there exists a unique Banach distance
 on $Y$, denoted $ dist_Y(y_1,y_2)=||y_1-y_2||_Y$ and  called  Federer-Whitney  $f$-Descendant  of the metric $dist_X$, such that the map $f$ is isometric at infinity with respect
 to $dist_X$ and $dist_Y$.}

\vspace {1mm} 

{\it Proof.}  Clearly, there exists a unique {\it maximal} Banach distance $dist_Y$ in $Y$  with respect to which $Lip_\infty(f)\leq 1$.

 To see that $coLip_\infty(f)=1$ as well, observe that every two disjoint $r$-balls in $Y$ satisfy
 $$dist_Y(B_{ y_1}(r), B_{ y_2}(r) )=\lim_{\lambda \to \infty} \lambda^{-1}\cdot dist_X 
\big (  f^{-1}(B_{ \lambda y_1}(\lambda r)),   f^{-1}(B_{ \lambda y_2}(\lambda r) )\big),$$
 where the scaling $y\to \lambda y$ is understood with some point in $Y$ taken for zero and 
 where the existence of the limit, as well the convexity and homogeneity of the limit distance function, trivially follow from

$\bullet$ the triangle inequality for $dist_X$,

$\bullet$ elementary properties of {\it asymptotically sub-additive} positive functions $d(\rho)$, $\rho>0$,  i.e.  such that  $d(\rho_1+\rho_2) \leq d(\rho_1)+d(\rho_2)+o(\rho_1)$;
the essential property of such functions $d$ is 
the existence of  the limit $\lim_{\rho\to \infty}d(\rho)/\rho$.  
 
 \vspace {1mm}

{\it Remarks.} (a) The above is the special (trivial) case of the full {\it  Federer-Whitney  duality  theorem}  between the {\it stable volume norm on homology $H_k(V)$} and the {\it  comass norm} on cohomology $H^k(V;\mathbb R)$.
(See sections 4.C and  4.D  in \cite {metr} and 7.4 in \cite {fil} for the statement and  further  applications of this theorem.)

(b)  The   Federer-Whitney theorem generalizes   \cite {pan} to maps $X\to Y$, where $Y$ is a simply  connected {\it nilpotent} Lie group with an {\it expanding} automorphism $\lambda$  and  where the existence of the limit  
$$\lim_{N \to \infty} \lambda^{-N}\cdot dist_X 
\big (  f^{-1}(\lambda^N(B_{ y_1})),   f^{-1}(\lambda^N(B_{  y_2} ))\big)$$
 is not so obvious.

 (c) The   property of being  isometric at infinity for equivariant maps  $f$ sometimes implies 
a much stronger (and non-obvious) one: 
$$|dist_X(f^{-1}(y_1), f^{-1}(y_1))-dist_Y(y_1,y_2)|\leq const <\infty.$$
This was proven by D. Burago \cite {bur} for surjective  maps from length metric spaces (where the distance is given by infima of length of curves between points)
into linear spaces $Y$ and, later on,  for maps into  some nilpotent Lie groups and 
 all hyperbolic $\Gamma$-spaces  by S. Krat  \cite {kra}, but  also there are counterexamples \cite {bre}. Yet,  the full geometry of  isometric at infinity equivariant maps  has not been fully clarified at the present day.

\vspace {1mm}
{\it  Abelian Volume Growth Inequality of  Burago-Ivanov .}  Let $\tilde V$ be  the universal covering of a compact  Riemannian $n$-manifold $V$ that admits a continuous map of non-zero degree onto the $n$-torus,  $V\to \mathbb T^n$,  e.g.  $V$ is  homeomorphic to   $ \mathbb T^n$.

{\it Then the  asymptotic volume growth of $R$-balls $B(R)$
 in  $\tilde V $ for $R\to \infty$ 
  is minorized by that in $\mathbb R^n$,}
   $$ \liminf_{R\to \infty} vol(B(R))/vol(B_{\mathbb R^n}(R))\geq 1.$$

{\it  Proof.} Let $Y$ be the universal covering of  the torus $\mathbb T^n$, let 
 $X$  be the corresponding $\mathbb Z^n$-covering of $V$ and $f: X\to Y$ the 
 $ \mathbb Z^n$-equivariant  map induced by our  $V\to  \mathbb T^n$.
   
 Let $||...||_Y$ be the Federer-Whitney norm on $Y$ and $h_\diamond$ John's quadratic Hilbertian (Euclidean) form. 
 
Since  $f$ is isometric at infinity with respect to  the Federer-Whitney metric $||y_1-y_2||_Y$,  the  $f$-pullbacks of the (Euclidean) $h_\diamond$-balls $B_\diamond(R-o(R))\subset \mathbb R^n =(Y, h_\diamond)$ of radii
$R-o(R^n)$ are contained in the Riemannian $R$-balls in $X=\tilde V $, and since
$Lip_\infty(f)=1$, the   asymptotic Banach-John volume inequality from the previous section applies.

  \vspace {1mm}
 
{\it  Remarks}. (a)  This proof:  

\vspace {1mm}

\hspace {7mm} Federer-Whitney + John's $h_\diamond$ + Hadamard's inequality
 
 \vspace {1mm}

\hspace {-6mm} is  similar to
 the original one  in \cite {bu-iv}, except that   
  the authors of  \cite {bu-iv} use at some point a rather subtle  Burago's theorem  (stated in the above  Remark (c))  instead of  the   (almost obvious) Federer-Whitney theorem.
   
 (b)    {\it Volume Rigidity.} The above argument combined with "sharpness of Hadamard"  (see \textbf {B} in section 2) implies by a simple argument  another theorem from \cite {bu-iv}.
 \vspace {1mm}

 Let  $X$ be a Riemannian $\Gamma$-manifold of dimension $n$, where the asymptotic volume growth of balls  equals that in $\mathbb R^n$,
$$ \liminf_{R\to \infty} vol(B_X(R))/vol(B_{\mathbb R^n}(R))=1.$$

{\it If $X$ admits a $\Gamma$-equivariant map   $f_0:X\to \mathbb R^n $ of a non-zero degree, where the action
 of $\Gamma$ on $\mathbb R^n$  is affine  co-compact quasi-translational, 
then $X$ is isometric to $\mathbb R^n$.}

 \vspace {1mm}

 This, however, leaves open the following  

{\it Question.} What is the minimal volume growth of an $X$ with the corresponding Federer-Whitney Banach space being in a given isometry class?

(We shall address this in Part 2 of this paper.) 
 
  \vspace {1mm}

 (c)  An attractive feature of the volume growth theorem,  is that   a general convexity  argument (John's e theorem)  yields, rather unexpectedly, a  {\it sharp  purely Riemannian} inequality. 
 
 (A non-sharp lower bound  on asymptotic volumes of balls follows by the above argument, see \cite {bab} and 4.C in \cite {metr}, from Federer-Whitney and  the inequality $|| ...||_A \geq \sqrt{n}\cdot ||...||_{h_\diamond}$ that  is the standard corollary to John's   
 $Hilb_{h_\diamond} \leq 1$.)
    
But  we shall see with our definition of  {\it Hilbertian volume} in section 9 that, in truth, there is not so much "Riemannian" in this inequality  after all: the Abelian volume growth inequality and the volume rigidity hold for (almost)  arbitrary metric $\Gamma$-spaces $X$.

\section {$\Gamma$-Spaces, Averaging and  Hilbert's Straightening Map $f_\diamond$.}

The  Federer-Whitney  Banach metric is, obviously, 
  invariant under actions of a group $\Gamma$ acting on our spaces.
Let us explain how to make   other  constructions  $\Gamma$-invariant as well. 

\vspace{1mm}

Recall, that  {\it a $\Gamma$-space}  in a given (topological) category is a space $X$ with an action of  a group $\Gamma$, denoted  $\gamma_X:X\to X$,  $\gamma \in \Gamma$, where these maps $\gamma_X$ are morphisms in our category.  We are mainly concerned with {\it metric $\Gamma$-spaces} where the maps  $\gamma_X$ are
isometries.

\vspace {1mm}

{\it Equivariant Hilbert Straightening Theorem.}  Let  $X$ and $Y$ be topological $\Gamma$-spaces as in the  Federer-Whitney theorem, i.e. $X$ is  a locally compact metric  $\Gamma$-space, and   $Y$  an $\mathbb R$-affine  $n$-space isomorphic to $\mathbb R^n$, where the action
 of $\Gamma$ is co-compact quasi-translational. 
 We denote by $Y_\bullet$ the linear space associated to $Y$ that is the space of parallel translations of $Y$ and that can be (non-functorially)  thought of  as $Y$ with a  point in it 
 distinguished for $0$.

 Let $f_0 : X \to Y$ be a proper   continuous $\Gamma$-equivariant map. 
 
\vspace {1mm}

{\it  Then there exist 
a unique   $\Gamma$-invariant  (Federer-Whitney-John's) Euclidean (Hilbertian)  metric  
$dist_\diamond=||...||_{h_\diamond}$ on $Y$, 
and  a (non-unique)  $\Gamma$-equivariant Lipschitz map $f_\diamond : X\to Y$,  such that
$$dist_{h_\diamond}(f_\diamond,f_0)\leq \infty, $$
   the Hilbert constant of $f_{\diamond}$ with respect to the Hilbertian form $h_{\diamond}$ on $Y$ satisfies
$$\widetilde {Hilb}(f_{\diamond})=1,$$}
and, at the same time,
$$coLip_\infty(f_{\diamond})\leq 1.$$}

{\it Proof.}  
The distance $||...||_{h_\diamond}$ in $Y$ for John's  $h_\diamond$ associated to  the Federer-Whitney  Banach metric  $dist_Y$ is  
$\Gamma$-invariant due to  uniqueness of $h_\diamond$  but the underlying measure $\mu_\diamond$ -- partition of unity ($d_{\mu_\diamond} p$) or equivalently,
partition of   $h_\diamond$ (that is $d_{\mu_\diamond} l$), is not necessarily invariant.

 But since the  linear isometry group $G$ of $(Y_\bullet, ||...||_Y) $, where $Y_\bullet$ is  the linear space  associated to $Y$, is compact, one can average partitions  of unity $\mu_\diamond$ (or rather partitons of   $  h_{\diamond}$ into squares of linear functions)   over $G$  and have John's partition of unity $\mu_\diamond$  invariant under the action of $\Gamma$.

 Let us  make the $\varepsilon$-straightened  maps $f_\varepsilon : X\to Y $  from section 3  equivariant  by the following   standard averaging over $\Gamma$.

Observe that  $\Gamma$ acts on the space $F$ of maps $f:X\to Y$ by 
$\gamma_F:f\mapsto \gamma^{-1}_Y\circ  f\circ\gamma_X$ with    fixed points of this action corresponding to equivariant maps $f$.

If we  replace maps  $f \in F$ by $b_f=f-f_0$ for our (equivariant!) $f_0:X\to Y$,
 this actions  becomes the obvious shift  action on the space of maps  $b: X\to Y_\bullet$ for $\gamma: b(x)\mapsto b(\gamma_X(x))$. 

Since the full isometry group $iso(Y)\supset \Gamma$ is {\it amenable}
one can average {\it  bounded}   functions to $\Gamma$-invariant ones,
where  "bounded  for $b_f$" corresponds to the "finite distance from $f_0$."
$$dist_Y(f,f_0)=\sup_{x\in X}dist_\diamond(f(x),f_0(x))<\infty.$$

Since the straightening maps  $f_\varepsilon: X\to Y$ from section 3 
 do lie within finite distance from $f_0$,  they can be averaged to equivaruant one, say $\overline f_\varepsilon : X\to Y$,  and since  $||\widetilde {dil}^\ast f_\varepsilon ||_{L_2(\mu)} \leq 1+\varepsilon$ and 
 the function(al)  $ f\mapsto ||\widetilde {dil}^\ast||_{L_2(\mu)}^2$ is (obviously) {\it convex},  these averaged maps also have their $\tilde L_2$-dilations   bounded by 
$$||\widetilde {dil}^\ast \overline f_\varepsilon||_{L_2(\mu)} \leq 1+\varepsilon$$

Finally, since the maps $\overline f_\varepsilon$ are uniformly Lipschitz, some sequence of them converges with $\varepsilon \to 0 $ to the   desired  equivariant  map $f_\diamond: X\to Y$. QED.

\vspace {1mm}

\vspace {1mm}

 {\it Remark}. The Hilbert straightening  map $f_\diamond$ can be regarded as a metric counterpart to the { \it Abel-Jacobi-Albanese maps} from K\"ahler manifolds to their Jacobians. The latter maps,
being holomorphic are (pluri)harmonic; they minimize  Dirichlet's $\int_{X/\Gamma} ||Df(x)||^2_{L_2}dx$
over all $\Gamma$-equivariant maps $f$,  while the map $f_\diamond$
minimizes, in a way,  the norm $\sup_x ||Df(x)||_{L_2} $.

Besides $f_\diamond$, there is also a {\it dynamical} counterpart to
 the Abel-Jacobi-Albanese, called  {\it the  Shub-Franks map} that is associated with {\it hyperbolic endomorphisms} of tori and infra-nil-manifolds in general.  The  three maps  and their generalizations intricately  intertwine in 
their applications to the geometric rigidity theory  (See \cite {shoe} 
for a related discussion on   A-J-A -- S-F connection.)

\subsection { On Infinite Dimensional $\Gamma$-Spaces and Codiffusion Spaces.}

The above suggests the following  $\Gamma$-version of the notion of the $L_q$-dilation from section 1.

Let $Y$ be  a Banach space  with  an affine isometric action of a locally compact (e.g. discrete) group $\Gamma$. Let  $\Pi\subset \Gamma$ be the subgroup of parallel translations in $\Gamma$ and let $\Gamma_0=\Gamma/\Pi$. 

Thus,   $\Gamma_0$ and, hence, $\Gamma$ itself,  act on $Y$ by {\it linear} isometric transformations.

{ \it A proper  axial  $\Gamma$-partition of unity in} $Y$ is a  {\it $\Gamma_0$-invariant} measure
 $ \mu=d_\mu p$ on the spaces $\cal P$ of axial projectors $p: Y\mapsto \mathbb R_p=im_p\subset Y$, such that the action of    $\Gamma_0$
 on the measure  space $({\cal P},\mu)$ is {\it proper}: 
 
 \vspace {1mm}
 
 there exists a measurable subset $U\subset {\cal P}$ such that the $\Gamma_0$-orbit
 of $U$ equals almost all ${\cal P}$ and such that the subset of those $\gamma\in \Gamma_0$ for which $\mu (\gamma(U)\cap U)>0$ is precompact in $\Gamma_0$.
 
  \vspace {1mm}
 
 If the action is proper, the quotient space
 $\underline {\cal P}=  {\cal P}/\Gamma_0$ carry a  Lebesgue-Rochlin measure, say $\underline \mu$, and one can integrate invariant functions on $\cal P$
 over  $(\underline {\cal P},\underline \mu). $

 Thus, one can define the  $L_q$-dilation of a $\Gamma$-equivariant map $f$
 from a metric $\Gamma$-space $X$ into $Y$ 
as 
 $$||dil^\ast f||_{L_q(\underline\mu=d_{\underline\mu}  p)}=\left (\int _{\underline{\cal P}}Lip^  q(f)d_{\underline \mu}  p\right)^{1/q}.$$

\vspace {1mm}

 Several  constructions we met earlier, namely  those which  do  not involve integration apart from that
over ${\cal P}$ generalize to this setting: that are

$\bullet$ Equivariant minimal  "norms" $||min.dil^\ast f||_{L_q}$ and  $||min.\widetilde {dil}^\ast f||_{L_q}$,

$\bullet$ Banach straightening at infinity
(section 3)  and equivariant Hilbert straightening, 

$\bullet$ John's $h_\diamond$, where one needs to(?), regretfully, assume that the Banach space $Y$ is "essentially Hilbertian"  to start with: it admits  a $\Gamma_0$-invariant Hilbert form $h$,
such that $||...||_Y\leq ||...||_h\leq const \cdot ||...||_Y$.

 (This assumption rules out, for example, the only natural candidate  for John's $h_\diamond$ on the spaces $l_q(\Gamma)$, that is Hilbert's  $h=||...||^2_{l_2(\Gamma)}$  for countable groups $\Gamma$.)

\vspace {1mm}

Where are we   to go from this point?
 Are there  unexplored  domains populated by  interestingly  structured infinite dimensional $\Gamma$-spaces ?   

One may start a search for them among  {\it symbolic $\Gamma$-spaces},
 \cite{topo}, \cite{manifolds},   and/or {\it concentrated spaces}  \cite {spaces}).  

\vspace {2mm}

{\it Axes, Non-Linear Partitions of Unity   and Codiffusion.}  A distance minimizing geodesic in a metric space $Y$ that is the image of an isometric embedding $\mathbb R   \mapsto \mathbb R_\bullet = im_p\subset Y$
admits a $1$-Lipschitz projection $p$ onto itself. Thus "many" spaces $Y$
admit "many" axial projectors $p: Y\to \mathbb R_p
\subset Y$. 

Notice that an axis $\mathbb R_\bullet \subset Y$ may have several 
projections $Y\to  \mathbb R_\bullet $. Two,  apparently, most natural ones
correspond to the two {\it Busemann functions} $b_\pm$ on  $Y$ associated to this axes with two possible orientations: $y\mapsto b_\pm(y)\in \mathbb R=\mathbb R_\bullet\mapsto Y$.

There are two alternative ways to define the $L_2$-dilation of a map $f:X\to Y$. 

\vspace {1mm}

(1) If $Y$ is a Riemannian/Hilbertian manifold with the Riemannian quadratic  differential form $g$,
 one may forfeit projectors  and use   measures $\mu=d_\mu b$ on the space $\cal B$ of $1$-Lipschitz functions $b: Y\to \mathbb R$.  

 The  normalization condition of  $\mu$ being a partition of unity may be formulated  by
requiring the integral of the squares of the differentials of our functions to be equal to $g$,
$$ \int_{\cal B} db^2d_\mu b= g.$$ 

Here it seems   reasonable  to limit the support of relevant measures $\mu$ to the space of Busemann functions or, more generally of all  {\it horofunctions}.  Besides, if $Y$ is a proper $\Gamma$-space, it may be to one's  advantage  to restrict to $\Gamma$ invariant measures $\mu $ on $\cal B$ and   integrate over  ${\cal B}/\Gamma$ rather than over $\cal B$.

\vspace {1mm}

(2) In order to sum/integrate projectors with a measure on the space of projectors one may use an additional structure in $X$,
called

 {\it Codiffusion.} Let  ${\cal N }={\cal N }(Y)$  be the space of   probability measures $\nu$  on $Y$ where $Y$ is embedded to ${\cal N }$  by assigning  Dirac's $\delta_y$ to all $y\in Y$.

A codiffusion  in $Y$ is a projection $R: {\cal N }\to Y \subset {\cal N }$, where
"projection" means $R^2=R$.

For example, every  complete convex (locally) affine topological space, e.g.  a closed convex subset in a Banach space, comes with a natural codiffusion that sends every measure $\nu$ on $X$
to its {\it center of mass}

The center of mass construction generalizes to those geodesic spaces, in particular to Riemannian manifolds, $Y$ where
the square distance functions $ y'\mapsto dist^2(y',y)$ are strictly convex
on all geodesics:   the center of mass of a $\nu$ is defined as the minimum point $y'_{min}\in Y$ of the (strictly convex!)
function $y'\to  \int_Y dist^2(y',y)d_\nu y$ (see \cite {gro-kar},  \cite {kar}  and also  \cite {afs} for the history of this concept and new applications.)

However, all this does not seem to help in answering  the following 

 {\it Questions}  (a) What are "interesting" instances of  equivariant maps $f: X\to Y$ between $n$-dimensional Riemannian (Finsler?)
 $\Gamma$-spaces such that $$coLip_\infty(f)\cdot \widetilde{Hilb}(f(x))=1\mbox { for all $x\in X$}$$
  where the group $\Gamma$ is {\it  not} virtually Abelian?
 
 (b) Let $X$  be a metric $\Gamma$-space,  $Y =(Y,h)$ a Finsler (e.g. Riemannian) $\Gamma$-manifold and $f: X\to Y$ be  a proper equivariant map. 
Suppose that all geodesics in $Y$ are distance minimizing. 
 
 What is the minimal $C= C(Y)>0$ for which 
 there exists a  $\Gamma$-equivariant map $f_\bullet:X\to Y$ within bounded distance from $f$ such that $\widetilde{Hilb}(f_\bullet(x)) \leq C\cdot Lip_\infty(f)$ for all
  $x\in X$?

 (c)  Let $Y$ be a contractible metric $\Gamma$-space. Consider all,
  possibly singular,  (oriented?) Riemannian  $\Gamma$-spaces $X$ of dimension $n$ that admit  equivariant maps $f:X\to Y$ that have $Lip_\infty(f)\leq 1$ and 
 such that the induced homomorphisms on their invariant $n$-dimensional real cohomologies do not vanish, e.g.   $dim(Y)=dim(X)=n$ and $deg(f)\neq 0$. 
  
  When does there exist an $X_{min}$ among them that
  minimizes $vol_n(X/\Gamma)$?

 How does $X_{min}$ change (if at all) with passing to subgroups of finite indices in $\Gamma$?
 
  Can this $X_{min}$ be explicitly described for
  particular "simple"  (possibly, infinite dimensional, e.g. isometric to Hilbert's $\mathbb R^\infty$) spaces  $Y$?

 Given  subgroups  $\Gamma_i\subset \Gamma$ of finite indices $ind_i=card(\Gamma/\Gamma_i)$,   let 
 $f_i:X_i\to Y$ be an  $ind_i^{-1}vol_n(X/\Gamma_i)$-minimizing  sequence of $\Gamma_i$-equivariant maps $f_i:X_i\to Y$ with $Lip_\infty(f_i)\leq 1$. When and how does such a sequence converge?

\section {Burago-Ivanov's Solution of Hopf Conjecture.}

Let us show (essentially)  following the original argument in \cite {bu-iv-h} how Hilbert's straightening implies {\it Hopf conjecture}:

 \vspace {1mm}
 
 {\it Let $V= (V\underline,  g)$ be a compact  Riemannian manifold without conjugate points
 and   $\underline f: V\mathbb \to \mathbb  T^n$ be a continous map such  that the induced homomorphism of the fundamental groups, $\pi_1(X)\to \pi_1(\mathbb T^n)=\mathbb Z^n$, is an isomorphism.   Then    the universal covering $X=(X, g)$ of $V$ is isometric to $\mathbb R^n$}.

\vspace {1mm}

 In fact,  if all geodesic segments in  $X$ are distance minimizing, then 
 
 \vspace {1mm}

{\it Hilbert's straightening map $f_\diamond:(X, g) \to (Y, h_\diamond)=\mathbb R^n$ is an isometry,}
 
 \vspace {1mm}

\hspace {-6mm} where $Y$ is the universal covering  of the torus $\mathbb T^n$  and  $h_\diamond$ is John's Hilbertian form  associated to the Federer-Whitney  $f$-descendant   $ ||...||_Y$ of the metric $dist_g$ on $X$ for $f :X\to Y$ being the lift of $\underline f$ to $X$.

\hspace {-6mm}  This is
  immediate with the following  simple classical formulas  \textbf {(LS1)} and \textbf {(LS1)}.

\vspace {1mm}

{\it Liouville-Santalo Integral Identities.} Let $X= (X,g)$  be a  complete  Riemannian  manifold.

Given a measure $\lambda =d_\lambda \tau$ on  the tangent bundle $T(X)$,  denote by  $\underline \lambda=d_{\underline \lambda}x$ the push-forward of $\lambda$ under the projection $\pi: T(X)\to X$, and assume  that $\underline \lambda(W)<\infty$ for all compact subsets $W\subset X$.

Call a measure $\lambda$ {\it balanced } if, for every continuous quadratic form $h$ on $T(X)$,
$$\int_X trace _{g_x}(h_x) d_{\underline \lambda} x=dim(X)\int_{T(X)} h(\tau) d_\lambda\tau. $$

For example, the  {\it Liouville  measure}, that is supported on the unit tangent bundle $S(X)\subset T(X)$, is balanced, since
the trace of a quadratic form $h$ on $\mathbb R^n$ equals $n$-times the average
 values of $h$ on the unit sphere $S^{n-1}\subset \mathbb R^n$.

Let  $ \Sigma(R)=\Sigma(X;R)$ be the space of geodesic segments of length $R$ in $X$ that are locally isometric (geodesic)  maps $\sigma: [0,R]\to X$ and denote by 
  $\partial_\sigma\in S_{\sigma(0)}(X) \subset T(X)$  the unit tangent vector to our geodesic $\sigma$
at the point $\sigma(0) \in X$.

 The correspondence  $ \sigma \leftrightarrow s=\partial_\sigma$  identifies   $ \Sigma(R)$ with  the unit tangent bundle $S(X)$; accordingly, the measures on  $\Sigma( R)$ corresponding to   $\lambda=d_\lambda s$ on $S(R)$ are denoted $d_\lambda \sigma$.

\vspace {1mm}

Let   $X$ be acted upon  discretely and  isometrically by a group $\Gamma$, 
and let $\lambda$ be a  measure on $S(X)$ that is {\it invariant  under the action of $\Gamma$} and that is {\it normal} in the sense
the total $\lambda$-mass of  $S(X)/\Gamma$ equals one,
  $$\lambda(S(X)/\Gamma)=\int_{S(X)/\Gamma} d_\lambda \tau=1.$$

 Let  $Y = (Y,h)$ be another  complete  Riemannian  $\Gamma$-manifold and   let $f :X\to Y$  be a proper   $\Gamma$-equivariant Lipschitz map.

Average/integrate the  dilation of $f$  at the two ends of our geodesic segments, or, rather,  the  dilation of the composed maps $f\circ \sigma: [0,R]\to Y$,
 $$dil_{f\circ \sigma}(0,R)=R^{-1}dist_Y(f\circ \sigma(0),f\circ \sigma(R)),$$
over $(\Sigma(R), d_\lambda\sigma)$ and denote this  averaged  dilation by
$$dil_{av}(f\circ \sigma,R)=  \int_{S(X)/\Gamma= \Sigma(R)/\Gamma}dil_{f\circ \sigma}(0,R)d_\lambda \sigma.$$

Clearly, this $dil_{av}$ is bounded by the {\it average length} of the parametized curves $f\circ\sigma :[0,R]\to Y$ divided by $R$, 
 $$dil_{av}(f\circ \sigma,R) \leq R^{-1}length_{av}(f\circ\sigma[0,R])= R^{-1} \int_ {\Sigma(R)/\Gamma} length (f\circ\sigma[0,R])d_\lambda \sigma.$$

If the measure $\lambda$ is {\it invariant under the geodesic flow}, then, by Fubini's theorem, this integral equals the integrated norm of  $\partial_\sigma f  
\in T_{f(x)}(Y)$, for $x=\sigma(0)=\pi (\partial_\sigma)$, that is the derivative (differential)  of $f$ in the direction of the unit tangent vector $\partial_\sigma \in S_x\subset T(X)$,

$$R^{-1}length_{av}(f\circ\sigma[0,R])=  \int_{\Sigma(R)/\Gamma} ||\partial_\sigma f || d_\lambda \sigma, \leqno { \textbf {(LS1)}} $$ 
where, by {\it Schwartz inequality,}
  $$ \int_{\Sigma(R)/\Gamma} ||\partial \sigma(0)|| d_\lambda \sigma \leq \left (\int_{\Sigma(R)/\Gamma} ||\partial \sigma(0)||^2 d_\lambda \sigma \right )^{1/2}.$$
If  $\lambda$ is {\it  balanced}, then 
$$\int_{\Sigma(R)/\Gamma} ||\partial \sigma(0)||^2 d_\lambda \sigma =Hilb^2_{av}(f),\leqno   { \textbf {(LS2)}} $$
where $Hilb_{av}(f)$ is the normalized {\it Dirichlet-Hilbert}  energy of $f$ that is
$$ Hilb^2_{av}(f)=_{def} n^{-1} \int_{X/\Gamma} ||Df||^2_{L_2}d_{\underline \lambda} x=n^{-1}\int_{X/\Gamma} trace_{g(x)} (Df(x))^\ast (h)d_{\underline \lambda}$$
for the pullback    $(Df(x))^\ast (h)$  of the Riemannian form $h$ from $T(Y)$ to $T(X)$ by the differential of $f$.
 
  \vspace {1mm}

In sum, { \textbf {(LS1)}}+{\textbf {(LS2)}}+(Schwartz inequality) imply:
  
  \vspace {1mm}
 
 {\it  if  $\lambda$   is a normal $\Gamma$-invariant  measure on $T(X)$ that is  balanced   and invariant under the geodesic flow and $f:X\to Y$ is a $\Gamma$-equivariant map, then the $\lambda$-average
  dilation of $f$ at the ends of geodesic segments in $X$ of length $R$ satisfies 
 $$ dil_{av}(f\circ \sigma,R) \leq Hilb_{av}(f) = n^{-1/2} \left(\int_{X/\Gamma} ||Df||^2_{L_2}d_{\underline \lambda} x\right)^{1/2}.$$}

 This integral inequality is sharp: the equality is possible only if  all quantities involved are equal point-wise and this sharpness is independent of $R$, since
 the integral $  \int_{\Sigma(R)/\Gamma} ||\partial_\sigma f || d_\lambda \sigma =R^{-1}length_{av}(f\circ\sigma[0,R])$ does not depend on $R$.
 
 It follows, in particular, that

 \vspace {1mm}

 {\it  if
  $$\limsup_{R\to \infty} dil_{av}(f\circ \sigma,R)\geq 1\mbox { and } Hilb_{av}(f)\leq 1,$$
 then the map $f$ is isometric on every geodesic $\sigma : (-\infty, +\infty)\to X$
 that is contained in the support of $\lambda$.}
 
 \vspace {1mm}

The Hopf conjecture follows  by applying the above
 to the universal covering $X$ of a torus without conjugate points with the Liouville measure $\lambda$ on $T(X)$ and to the  Hilbert straightening map $f_\diamond : X\to \mathbb R^n=(Y,h_\diamond)$.

In fact,  "no conjugate points" says that   all geodesic segments  in $X$, that are $\sigma:[0,R]\to X$,
have $dist_X(\sigma(0),\sigma(R))=R$. Therefore, 
the  average dilation of $f_\diamond$  satisfies
$$\big (\limsup_{R\to \infty} dil_{av}(f_\diamond\circ \sigma,R)\big )^{-1}\leq coLip_\infty(f_\diamond),$$
where, recall,
$$coLip_\infty(f)= \liminf_{dist(x_1,x_2)\to \infty} dist(f(x_1), f(x_2))/dist(x_1,x_2)$$
for proper maps $f$.

On the other hand the inequality $Hilb(f_\diamond) \leq 1$ implies that
$||Df_\diamond||^2_{L_2}\leq n$ for almost all $x\in X$.
Hence, $f_\diamond$ is an isometry. QED.

\vspace {2mm}
{\it Questions.}  Let $X$ be a complete  Riemannian $\Gamma$-manifold 
  of dimension $n$ and $f: X\to \mathbb R^n$ a proper continuous  
  $\Gamma$-equivariant map for some discrete cocompact isometric action of the group $\Gamma$ on $\mathbb R^n$. For example, $X$ may be equal the universal covering of a compact manifold homeomorphic to the   $n$-torus.

 What is a possible geometry of $X$ if
$$dil_{av}(f\circ \sigma,R) \geq Hilb_{av}(f) -\varepsilon,  \mbox { for all }  R\geq R_0,$$
where the averages  are taken with the Liouville measure?

What can be geometry of an  $X$ if  the geodesic segments
 $\sigma :[0,R] \to X$,  satisfy
$dist_X(\sigma(0),\sigma(R))\geq \phi(R)\cdot R$ for a given  function
$0<\phi (R)<1$ and all sufficiently large $R$?

In particular, does the inequality  $dist_X(\sigma(0),\sigma(R))\geq (1-\varepsilon)R$
for a small (any$<$1?) $\varepsilon$ and all $\sigma$
imply that $X$ is isometric to $\mathbb R^n$?

Do  "majority" of  $X$  have $dist_X(\sigma(0),\sigma(R))\sim R^{1/2}$ on the average?

\section { Hilbert Jacobian,   Homological Volume, Semicontinuity,  Coarea Inequality and John's  Identity.}   
 
 Let $V=(V^n,\chi)$ be a smooth $n$-dimensional manifold with a smooth measure $\chi$. Assume that $V$ is diffeomorphic to an open subset in $\mathbb R^n$   and, given a Lipschitz map $f$ from a metric space $X$ to $Y$,  consider all diffeomorphisms $\phi: V\to \mathbb R^n$ for which $Jac^{[n]}(\phi(v))\geq 1$ for all $v\in V$ (such 
$ \phi(v)$ always  exist and can be even taken smooth).  

Define
  the {\it global  Hilbert Jacobian} of $f$ restricted to a subset 
 $U\subset X$ as

$$\widetilde{Hil}Jac_{glb}^{[n]} (f|_U)=\inf_{\phi} \widetilde{Hilb}^n(\phi\circ f|_U), \mbox { }  \phi\circ f|_U:U\to \mathbb R^n$$

and 
$$\widetilde{Hil}Jac^{[n]}f_{/v}=\inf_{V_v\ni v} \widetilde{Hil}Jac_{glb}^{[n]} (f|_{f^{-1}(V)}) ,\mbox  {  }
y\in Y,$$
where the infimum is taken over all neighbourhoods $V_y\subset V$ of $v.$

Since every smooth $n$-manifold is locally diffeomorphic $\mathbb R^n$, this 
definition makes sense for {\it all } smooth, including closed ones,  manifolds  $V=(V^n, \chi)$.

Next let
$$\widetilde{Hil}Jac^{[n]}f(x)=\inf_{U\ni x} \widetilde{Hil}Jac_{glb}^{[n]} (f|_U) ,\mbox  {  } x\in X,$$
 where the infimum is taken over all neighbourhoods $U\subset X$ of $x$ and where, again, this is defined for all $V$.

 Finally, define the  point-wise  Hilbert Jacobian $\widetilde{Hil}Jac^{[n]}f(x)$ of  a Lipschitz map $f$   between arbitrary metric spaces, $f:X\to Y$,   as the infimum   of the numbers $J$
 such that 
 $$ \widetilde{Hil}Jac^{[n]}(\tau \circ f)(x)\leq J\cdot \widetilde{Hil}Jac^{[n]}f(x)$$
for all Lipschitz maps $\tau: Y \to \mathbb R^n$.
 
 It is immediate that 
 $$ \widetilde{Hil}Jac^{[n]}(f)\leq Lip^n(f)$$
 and that 
 composed maps $X\underset {f}\to Y \underset {g}\to Z$ satisfy
$$ \widetilde{ Hil}Jac^{[n]}g\circ f(x) \leq \widetilde{Hil} Jac^{[n]} f(x)\cdot  \widetilde{Hil}Jac^{[n]} g(y)\mbox  { for }y=f(x)$$

Also observe that  $\widetilde{Hil}Jac^{[n]}(f)$ equals the (absolure value of the)  ordinary  Jacobian $Jac^{[n]}(f)$ for smoth maps $f$ between Riemannian manifolds.

\vspace {1mm}

{\it  Homologically Stable Hilbertian Volume.}  We want to define $n$-volumes in metric spaces  
$X$ with a usual behaviour under (Jacobians  of)
Lipschitz maps $f:X\to Y=(Y^n,\chi) $,  e.g. decreasing under  maps $f$ with  $\widetilde{Hil}Jac^{[n]}(f(x))\leq 1$, where we are concerned with {\it stability of the images} under small continuous perturbations maps that is defined as follows.
 
 Fix a coefficient group (e.g. a field) $\mathbb F$
and call  a point $y\in Y$ {\it cohomologically ($H^n_\mathbb F$-image) stable} for a continuous map $f:X \to Y$ if for all sufficiently small neighbourhoods $V=V_y\subset Y$ of $y$  the induced  homomorphism  on the relative cohomology
$$f^{\ast n} : H^n(V,\partial V;\mathbb F)\to  H^n(f^{-1}(V), f^{-1}(\partial V);\mathbb F)$$
{\it does not vanish}.
The set of these stable points is denoted by $stbl_\mathbb F im_f(X) \subset Y$.

As we shall be dealing with locally compact spaces, the cohomology will be {\it  \v{C}ech;} to be specific, we stick to  $\mathbb F=\mathbb F_2=\mathbb Z/2\mathbb Z$. 

\vspace {1mm}

{\it The   {hyper-Hilbertian} global  $n$-volume   $glb_\mathbb F vol^{[n]}(
\mathbb R^n\backslash X)$ of $X$ over $\mathbb R^n$} is  the supremum of the numbers $H$ such that $X$ admits a Lipschitz map  $f: X\to \mathbb R^n$ with 
$\widetilde{Hil} Jac^{[n]}_x( f) \leq 1$, $x\in X$, and such that the Lebesgue-Haar measure of the set $stbl_\mathbb Fim_f(X) \subset \mathbb R^n$ is $\geq H$.

For example,  $glb_\mathbb F vol^{[n]}(
\mathbb R^n\backslash X)=vol_n(X), $ for smooth  open Riemannian manifolds  $X$,
 while closed manifolds $X$ satisfy
 $glb_\mathbb F vol^{[n]}(
\mathbb R^n\backslash X)=vol_n(X)/2. $
(To get rid of this silly $1/2$ one has to  map $X$ to the $n$-spheres  $S^n $ rather than to $\mathbb R^n$.)

A more interesting example is provided by closed minimal hypersurfaces $X=X^n$ in  compact Riemannian manifolds $X^{n+1}\supset X$ representing homology classes in $H_n(X^{n+1})$: such an $X$ admits a $\lambda$-Lipschitz map, $f: X\to S^n$ of non-zero degree, where $\lambda\leq \lambda_0(X^{n-1}, vol_n(X))<\infty$ as it follows from the {\it compactness} of the space of minimal subvarieties in $X^{n+1}$ with bounded volumes  and the {\it maximum principle} for minimal hypersurfaces.
 
Then one sees,  by looking at tangent cones, that every point $x$ in such a minimal $ X \subset X^{n+1}$ admits an arbitrarily small  neighbourhood $U_x\subset X$ such that  
 $$glb_\mathbb F vol^{[n]}(
\mathbb R^n\backslash U_x) \geq \varepsilon \cdot vol_n(U_x)\mbox
{ for some }  \varepsilon\geq\varepsilon (X)>0.$$

\vspace {1mm}

{\it Homology Domination and Semicontinuity of $glb_\mathbb F vol^{[n]}$.} Let $Z$ be a metric space,  $X\subset Z$ a locally compact locally closed (open $\cap$ closed) subset, and let $X_i$, i=1,2,..., be a sequence of locally closed subsets that converges to $X$ in the sense that for every neighbourhood $Z_X\subset Z$ of $X$
all but finitely many of $X_i$ are contained in $Z_X$.

Say that $X_i$ {\it (co)homologically dominate $X$ in dimension $n$ with a given  coefficient domain $\mathbb F$}
if for every pair of open subsets $Z_1\subset Z_2\subset Z$
every relative cohomology class in $H^n(Z_2,Z_1;\mathbb F)$  that
restricts to a {\it non-zero} class in $H^n(Z_2\cap X,Z_1\cap X;\mathbb F)$
also restricts to {\it non-zero} classes in  $H^n(Z_2\cap X_i ,Z_1\cap X_i;\mathbb F)$ {\it for all but finitely many} $i$.

By invoking the  local $\tilde L_q$-extension with $q=2$ from section 1, one obtains: 

\vspace {1mm}
{\it If  subsets $X_i \subset Z$ homologically dominate a subset $X$ in a metric space $ Z$ then
$$\liminf_{i\to \infty}glb_\mathbb F vol^{[n]}(X_i)\geq glb_\mathbb F vol^{[n]}(X).$$}

\vspace {1mm}

{\it (Semi)Local Coarea Inequality.} Let  $Y=Y^m $,  $m\leq n$, be a topological manifold and  $f:X\to Y$  a continuous map. Call the map $f$  {\it   (co)mologically (co)stable over a point $y$ in dimension $n$ with coefficients $\mathbb F$} if the   fiber $ X_y=f^{-1}(y)\subset X$ and $f$ near $X_y$ have the following property. 

 \vspace {1mm}
 
 \textbf {[STB]} Let $f_y:  X\to \mathbb R^{n-m}$ be a continuous  map and $V\subset \mathbb  R^{n-m}$ be an open subset such that the induced cohomology
 homomorphism 
 $$ H^{n-m}(V, \mathbb \partial V ;\mathbb F) \to H^{n-m}(f_y^{-1}(V)\cap X_y, f_y^{-1}( \partial V)\cap X_y;\mathbb F )$$
  does not vanish. Then for all sufficiently small closed  neighbourhoods 
  $W_y \subset Y$ of $y$ the map 
  $$F=(f,f_y): f^{-1}(W_y)\to W_y\times \mathbb R^{n-m}$$
  induces {\it non-zero} homomorphism 
 $$ F^\ast: H^n(W_y\times V, \partial( W_y\times V) ;\mathbb F)\to 
 H^n(F^{-1}(W_y\times V), F^{-1}(\partial( W_y\times V) ;\mathbb F)).$$
 
 For example an $\mathbb R$-valued  function $f$ on a topological $n$ manifold $X$ is stable over a point $y\in \mathbb R$ if the level $f^{-1}(y)$ contains  no local maxima and/or minima points of $f$.

  \vspace {1mm}
 
 Now let $Y$ be smooth $m$-manifold with a smooth measure $\chi^{[m]}$ on it,  let   $f:X\to Y$ be a Lipschitz map and let $W_y(\varepsilon) \subset Y$ be $\varepsilon$-balls around $y$ for some metric in $Y$.

{\it If $f$ satisfies} \textbf {[STB]}  {\it then  
$$\liminf_{\varepsilon \to 0} glb_\mathbb F vol^{[n]}(f^{-1}(W_y(\varepsilon))  \geq \frac{  \chi^{[m]}(W_y(\varepsilon))\cdot glb_\mathbb F vol^{[n-m]}(f^{-1}(y))  }{\widetilde{Hil}Jac^{[m]}f_{/y}}.$$}

{\it Proof}. The condition   \textbf {[STB]}, albeit  unpleasantly restrictive,   matches
  the  local $\tilde L_q$-extension property and the coarea inequality trivially follows.

\vspace {2mm}

Let us partially localize the above "volume" with countable systems  $\cal U$ of  open subsets $U_i\in X$, $i\in I$ that have bounded  intersection multiplicity, denoted $mult({\cal U}) <\infty$,
and let 
$$Hil.vol_\mathbb F^{[n]}(\mathbb R^n\backslash X)=\sup _{\cal U}
(mult({\cal U}))^{-1} \sum_{i\in I} glb_\mathbb F vol^{[n]}(
\mathbb R^n\backslash U_i),$$
where this  will be often abbreviated  as follows
 $$\widetilde{Hil}.vol^{[n]}(X)\mbox { instead of  } \widetilde{Hil}.vol_\mathbb F^{[n]}( \mathbb R^n\backslash X)$$
with  $\mathbb F=\mathbb F_2$, unless otherwise indicated.

The homological volume  $\widetilde{Hil}.vol^{[n]}(X)$ of smooth and piece-wise smooth Riemannian spaces $X$ equals the ordinary  volume $vol_n(X)$.  

But, in general,  $ U\mapsto \widetilde{Hil}.vol^{[n]}(U)$, $U\subset X$,  is not,  a priori,   an additive set function. This can be amended
but using countable systems ${\cal U}_j$ of open sets  $U_{ij}\subset X$ where $diam(U_{ij})\to 0$ for $j\to \infty$ and localizing with
$$loc.vol(X)= \limsup _{j\to \infty}{\cal U}_j
(mult({\cal U}_j))^{-1} \sum_{i\in I} glb_\mathbb F vol^{[n]}(
\mathbb R^n\backslash U_{ij}).$$
 
 This "local volume" equals $\widetilde{Hil}.vol^{[n]}(X)$ if   small $\varepsilon$-balls in $X$ have volumes $\geq V_\varepsilon$ 
 where $  V_\varepsilon/\varepsilon ^{n+1} \underset  {\varepsilon \to 0}\to \infty$;
 in general, 
   however, the local volume, probably,   may be (I have not worked out a convincing example)  smaller than   $\widetilde{Hil}.vol^{[n]}$  that is   defined  with unrestricted $U_i$.    To compensate  for this, we shall introduce another "volume"  in Part 2 of the paper   by means of maps $X\to \mathbb R^\infty$ with controlled filling volumes of the boundaries of images of subsets (or rather of $n$-cycles) similarly in certain respects  to \cite {amb}.

\vspace {1mm}

\vspace {1mm}

{\it  Jonn's-Hilbert Volume of Banach Balls.} The obvious corollary to John's Ellipsoid Theorem (pointed out in slightly different terms in \cite {bu-Iv})
reads: 

\vspace {1mm}

{\it the Hilbertian volume of the unit ball $B_X(1)$ in an $n$-dimensional Banach space $X$ equals its Euclidean volume with respect to John's Euclidean (Hilbertian) metric $h_\diamond$,}
$$Hil.vol^{[n]}(B_X(1)) =vol_{h_\diamond}(B_X(1)),$$
Consequently, since $B_{h_\diamond}(1) \subset B_X(1)$,

\vspace {1mm}
{\it the Hilbertian volume $vol_{h_\diamond}(B_X(1))$ is  greater
 than the (ordinary) volume of the unit Euclidean ball,  where the inequality is strict unless
$X$ is isometric to 
$\mathbb R^n$}.


\section {Lower Volume Bounds.}

Let us reformulate and reprove  a few standard geometric inequalities in the Hilbert volume   setting where  we abbreviate: $vol^{[n]}(X)=  Hil.vol^{[n]}(X)$ for general metric spaces $X$.

\vspace {1mm}

{\it $\square$-Spaces and their Faces.} A topological space $P$  with a given 
collection of  $n$ pairs of disjoint closed subsets,  
$ \pm Q_i =\pm Q_i(P) \subset P$, $i=1,..,n$, is called  an {\it $n$-cubical space} where $ \pm Q_i$ are called {\it faces} and their union $\partial P=\bigcup_i\pm Q_i$ is regarded as {\it the boundary} of $P$.

 Every such $P$, assuming it is a normal space,  admits a 
continuous map $\pi $  into the $n$-cube, $\pi: P\to \square= [0,1]^n$, such that the faces  $\pm Q_i(P)$ go to the $(n-1)$-faces  $\pm Q_i(\square)$ of the cube, where such a map is given by $n$-functions $\pi_i: P\to [0,1]$ where $-Q_i$ goes to $0$ and $+Q_i$ goes to $1$.  Clearly, such a map  $\pi$ is unique up-to homotopy.

  \vspace {1mm}
  
  { \it $\square$-Straightening Lemma.} Let  $P$ be an $n$-cubical metric space,  
  and $\pi : P\to \square=[0,1]^n$ be a continuous map that sends faces of $P$ to the corresponding faces of the cube $ \square$.

{\it If
$$dist _P(Q_i,-Q_i)\geq  \lambda_i, \mbox { } i=1,..,n,$$
then 
  $\pi  $ is homotopic to a Lipschitz map $F: P\to \square= [0,1]^n$ given by $n$-functions $f _1,...,f_n:P\to [0,1]$, such that
$$ Lip(f_i) \leq \lambda_i^{-1}.$$}

{\it Proof}. Take $f_i(x)=\lambda^{-1}dist (x, Q_i)$ for all $x\in P$,where $dist (x, Q_i)\leq \lambda_i$
and let $f_i(x)=1$ for  $dist (x, Q_i)\leq \lambda_i$.  
  
  \vspace {1mm} 
   
   {\it  Corollary: Besikovich-Derrick-Almgren   $\square$-Inequality}. If the map $\pi : P\to \square=[0,1]^n$ has non-zero degree, i.e. the induced cohomology  homomorphism $\pi^\ast: H^n(\square, \partial \square;\mathbb F)\to H^n(P,  \partial P;\mathbb F)$ doe not vanish for some coefficient domain $\mathbb F$,  then

   {\it  the $\mathbb F$-homological Hilbertian volume of $P$ is bounded from below by the product of the distances 
   between opposite faces in $P$,
  $$vol^{[n]}(P) \geq \prod_{i=1,...,n}dist_P(Q_i,-Q_i).$$ } 

{\it Proof.} Since $F$ is homotopic to $\pi$  and   $deg(\pi)\neq 0$,
 the homomorphism $F^\ast: H^n(\square, \partial \square;\mathbb F)\to H^n(P, \partial P;\mathbb F)$ does not vanish; hence, all points $y\in Y$ are cohomologically stable and 
$$1=vol_n(\square)\leq   ||Jac^{[n]}(F)||_{sup}\cdot vol^{[n]}(P) \leq vol^{[n]}(P)\cdot \prod_i \lambda_i^{-1}.$$  

{\it Example.}  Let  $X$ a Riemannian $n$-manifold and $\pm V_i  \subset X$, $i=1,2,...,n$, be smooth closed domains that are  $n$-submanifolds with smooth boundaries $\partial (\pm V_i) \subset V_i$, e.g. $X=\mathbb R^n$ and  $\pm V_i 
\subset \mathbb R^n$ are half-spaces.

 Let    $P $  equal the intersection of $\pm V_i $,
 $$P=\bigcap_i\pm V_i\mbox { and } \pm Q_i=\pm Q_i(P)=P\cap \partial (\pm V_i).$$

\vspace {1mm}
If $P$ is compact,  if  the $n$ hypersurfaces  $\partial( + V_i)$, $i=1,2,...,n$,  intersect {\it transversally}  and if the  intersection 
$$\bigcap_i+ Q_i= P\cap  \partial( + V_1)\cap \partial( + V_2)\cap...\cap \partial( + V_n)$$
(that is necessarily finite)  consists of  {\it odd} number of points, then 
$$vol_n(P) \geq \prod_i dist(+Q_i, -Q_i).$$

Indeed, our  map $F:  P\to \square$ has non-zero $\mathbb F_2$-degree because
the $F$-pullback of every  point $y\in \square$ close to the corner $\bigcap_i+Q_i(\square)$ of  $\square$ is  a finite set of  odd
cardinality.

\vspace {1mm}

Recall that  the map $\pi: P\to \square$, for every $n$-cubical space $P$ sends $\partial P\to \partial \square$ and
the essential property of $\pi$ used  in the construction of $F$, say for 
$dist_P(-Q_i.+Q_i)\geq 1$, is that $\pi$ is {\it $1$-Lipschitz on the boundary} $\partial P$
with respect to the {\it $sup$-norm} on $\mathbb R^n\supset [0,1]^n=\square$,
that is $||(y_1,...,y_i,...,y_n)||_{sup}=max_i|y_i|$.

Seen from this angle, the $\square$-inequality appears as a special case of the following corollary (stated slightly differently in  \cite {bu-iv})  of  Jonn's  theorem.   
\vspace {1mm}

{\it Filling Extremality of Banach spaces.}  Let $A=(A||...||_A)$ be an $n$-dimensional Banach space  and  let $U\subset A$ be a bounded open subset with connected boundary $\partial U$.

 Let $X$ be a  compact metric space, let
$Y\subset X$  be a closed subset and let  $F: Y\to \partial U\subset  A$  
be a $1$-Lipschitz map.

\vspace {1mm}

{\it If the inclusion/restriction homomorphism  $H^{n-1}(X;\mathbb F)\to H^{n-1}(Y;\mathbb F)$ vanishes  for some coefficient group $\mathbb F$, while   the cohomology homomorphism $ f^\ast: H^{n-1}(\partial U;\mathbb F)\to H^{n-1}(Y;\mathbb F)$
does not vanish,  then
$$vol^{[n]}(X) \geq vol^{[n]}(U)=vol_{h_\diamond}(U)$$
for  John's Euclidean (Hilbertian)  quadratic form 
$h_\diamond$ on $A$.}

\vspace {1mm}

{\it Proof}. Extend $f$ to a Lipschitz map $F:X\to A$ with $Hilb_{h_\diamond}(F)\leq 1$  (see   section 1) and observe that
the cohomology assumptions imply that this map is homologically image stable at all points $u\in U$.

Since 
 the Hilbertian Jacobian of $F: X\to (A, h_\diamond)$ is bounded by 
$$Hil.Jac^{[n]}(F) \leq Hilb^n(F)\leq 1,$$
the proof follows from John's volume identity in section 8.

 \vspace{2mm}

{\it $\triangle$-Inequalities.} There is  another generalization of the $\square$-nequality, besides the above  Banach  filling extremality, where the cube $[0,1]^n$ is replaced
by the Cartesian product of regular  Euclidean $n_j$-simplices,  denoted $\times _j \triangle_{n_j}$, $j=1,...,k,$.

Start with the case $k=1$ and define a {\it $ \triangle_{n}$-space} as 
 a topological space $P$ with a distinguished collection of subsets, called { \it faces} 
$Q_i=Q_i(P)$, $i=0,1,...n$, indexed by the $n+1$ codimension one faces in $\triangle_{n}$, such that the intersection $\bigcap_iQ_i$ is empty and where the union $\bigcup_i Q_i \subset P$ is regarded as the boundary $\partial P\subset P$

If $P$ is a normal space, it admits a continuous map $\pi: P\to\triangle_{n}$, such that  every face of $P$ goes to the corresponding face of the simplex
 $ \triangle_{n}$, where this $\pi$ is unique up to homotopy.

 Let $P$ me a metric space and denote by $ \Sigma= \Sigma_P:
P\to \mathbb R_+$ the sum of the distance functions to the faces,
  $$\Sigma(p)=\sum_i dist_P(p,Q_i),$$
 let 
 $$\Sigma_{\partial}= \Sigma_{\partial}(P)=\inf_{p\in \partial P}\Sigma(p),$$
and let 
$$S=\Sigma_{\partial}(P)/ \Sigma_{\partial}(\triangle_{n}).$$
\vspace {1mm}

{\it If the map $\pi :P\to  \triangle_{n}$ has non-zero degree, then
  $$Hil.vol^{[n]}(P)\geq
   S^nvol_n( \triangle_{n}).$$}

{\it Proof.}  There obviously exits an  axial partition of unity  ${\cal P}=\{p_i\}$ $i=0,1,...,n$ in $\mathbb R^n \supset  \triangle_{n}$
with $n+1$ axes $\mathbb R_i \subset \mathbb R^n$ that are normal to the $(n-1)$-faces of  $\triangle_{n}$ and such that the quaisi-projectores $p_i$ on  
$\triangle_{n}$, that are
$p_i:\triangle_{n}:\mathbb R_i=\mathbb R$, are given by  the functions
$y\mapsto \lambda\cdot dist(y)$ for  $\lambda=\sqrt{n/(n+1)}.$

Let $f_i: P\to \mathbb R=\mathbb R_i$ be given by $ p_i\mapsto \lambda S^{-1}\cdot dist(p,Q_i)$ and  observe that the map
$$F=\sum_i f_i: P\to \mathbb R^n\supset \triangle_{n}$$ 
sends the boundary of $P$ {\it outside the interior of $ \triangle_{n}$} where, moreover, this map 
$$F:\partial P\to \mathbb R^n\setminus int (\triangle_{n})$$
is {\it homotopic to} 
$$\pi: \partial P\to \partial \triangle_{n} \subset  \mathbb R^n\setminus int (\triangle_{n}).$$

Since $deg(F)=deg(\pi)\neq 0$,  the image of $F$  homologically stably covers  $\triangle_{n}$ and the proof follows for
 $$Hil.Jac^{[n]}(F) \leq Hilb^n(F)\leq S^{-n}.$$

 \vspace {1mm}

Now we turn to metric  {\it $\times _j \triangle_{n_j}$-spaces} $P$, where $j=1,2,...,k$,  and where  a $\times _j \triangle_{n_j}$-structure is given by a collection of "faces" $Q_{ji} \subset P$, $i=i_j=0,...n_j$ and where $Q_{ji}$ correspond to 
those  faces of codimension $1$   in the Cartesian product that equal the pullbacks
of the faces of $\triangle_{n_j}$ under the  coordinate projection  $\times _j \triangle_{n_j}\to  \triangle_{n_j}$.

 The condition we impose on the set of faces   $Q_{ji} $ reads:  a subset of faces has a common point in $P$ only if the corresponding faces in $\times _j \triangle_{n_j}$ have a common point.
 
This is equivalent to the existence of a continuous map $\pi: P\to  \times _j \triangle_{n_j}$ that sends faces to faces, where such a $\pi$ is unique up-to a homotopy.

Let  
$$\Sigma_j(p)=\sum_{i=1,...,n_j} dist_P(p,Q_{ji}), \mbox { } 
\Sigma_{\partial_j}= \Sigma_{\partial_j}(P)=\inf_{p\in \partial P}\Sigma_j(p),$$
$$S_j=\Sigma_{\partial_j}(P)/ \Sigma_{\partial}(\triangle_{n_j}).$$

{\it If the map $\pi :P\to  \triangle_{n}$ has non-zero degree, then
  $$Hil.vol^{[n_1+...+n_j]}(P)\geq
  \prod_{j=1,...,k} S^{n_j}_jvol_{n_j}( \triangle_{n_j}).$$}
 
 {\it Proof.} Take a partition of unity in $\mathbb R^{n_1+...+n_k} \supset \times _j \triangle_{n_j}$  with the axes normal to the faces of  $\times _j \triangle_{n_j}$and argue as   earlier.

  \subsection{ Digression: Acute Polyhedra.}

  What should be a generalization of the above
  for {\it non-regular} simplices $\triangle_{n_j}$? 
   Namely, what is the sharp lower bound on the distance function $dist_P(p_1,p_2)$ for $p_1,p_2\in \partial P$ that would imply $Hil.vol^{[n]}(P) \geq\times_j \triangle_{n_j}$?

Notice in this respect, if $P$ is a  convex    {\it acute}  Euclidean polyhedron, i.e.  all dihedral angles are  $\leq \pi/2$ then the distances to its codimension $1$ faces 
 $d_i:p\mapsto dist(p, Q_i( P ))$ satisfy
 $$\sum_i c_id_i =const=vol( P )\mbox { for } c_i   = n^{-1}vol_{n-1}(Q_i(P))$$
 but the mere  $\sum_i c_id_i \geq const$ does not  non-sufficient  even for
 $P$  (homologically) over (i.e.  mapped  with non-zero degree to) {\it non-regular} Euclidean triangles $\triangle_2$.

What is the  minimal  set of (preferably linear and/or log-linear)  inequalities between $d_i$ that would imply  $Hil.vol^{[n]}(P) \geq\times_j \triangle_{n_j}$?
 
This  question make sense for all acute Euclidean polyhedra besides 
$\times_j \triangle_{n_j}$ (some "non-acute" lower volume bounds can be derived from "acute" ones , see  7.3 in  \cite   {fil} ) but, in fact, there aren't any other by  {\it Steinitz}(?) {\it theorem}:

 \vspace {1mm}
 
 {\it every acute  polyhedron $P\subset \mathbb R^n$  is a Cartesian product of simplices.}

\vspace {1mm}
 {\it Proof.} Start by observing   that an  acute convex {\it spherical polyhedron} $P$  (that is an   intersection of hemispheres in $S^n$)
 is a  either a {\it simplex}, or, in the degenerate case, the {\it spherical suspension} over a simplex
in $S^{n-i}\subset S^n$.

 Indeed,  the dual polyhedron, say $P^\perp_p \subset S^{n}$, 
    has all its edges longer than $\pi/2$.  Consequently, the distance between 
    {\it every two} vertices in  $P^\perp_p$ is $\geq \pi/2$; hence, there are  at most $n+1$ vertices
     in  $P^\perp_p$.

Also note  that (non-strictly)  acute spherical triangles $\triangle \subset S^2$ have all their edges
bounded in length by $\pi/2$. It follows that all $m$-faces, $m=2,3,..., n-1$, of acute spherical
 $n$-simplices are acute.

\vspace {1mm}

 Now,  let  $P$ be {\it Euclidean} acute      
  take an $(n-1)$-face $Q \subset P$, 
   move its supporting hyperplane $ H^{n-1}=H_Q \supset Q$, 
    inward parallel to itself  until it hits a vertex, 
    say $p\in P$ and denote by  $  H_p\subset \mathbb R^n$
    the so moved hyperplane.    
  
  Since the normal projections from $P$ onto the hyperplanes supporting  $(n-1)$-faces adjacent to $Q$ send $P$ {\it into} (hence, onto) these faces, 
  all of $P$ is contained in the "band"  $[H_Q, H_p] \subset \mathbb R^n
        $ between  the parallel hyperplanes $H_Q \subset \mathbb R^n $ and $H_p \subset \mathbb R^n$.        
       Moreover,   
     
     $\bullet $ if  the opposite face $-Q=H_p\cap P$ has $dim(-Q)=n-1$ then  $P$ is orthogonally splits into the Cartesian product of 
     $Q$ with a segment; it is seen with the two orthogonal projections $P\to Q, -Q$.
 
 $\bullet \bullet $ If $Q$ is an $(n-1)$-simplex and $dim(-Q)<n-1$ then, obviously, $P$ is an 
 $n$-simplex.

      Since the spherical 
 polyhedra underlying the tangent cones at all vertices
     of $P$ are acute,  $\bullet $  and $\bullet \bullet $    apply  to  the  faces  of  $P$.  
         
  Namely, let $\Delta^m \subset P$ be a {\it simplex}-face  of {\it maximal} dimension $m$.       
   Then every $(m+1)$-face containing this    $\Delta^m$   orthogonally splits according to $\bullet$ and $\bullet \bullet$.  It follows  by induction on $m$ that all of $P$ splits, 
   $P=\Delta^m\times P'$, and the proof is concluded
   by induction on $n$.
       
      \vspace {1mm}

Conclude by noticing  that the number of $(n-1)$-faces of a  {\it convex} polyhedron  
$P^n\subset \mathbb R^n$  with all dihedral angles $\leq \pi-\varepsilon$,  is (obviously) universally bounded, say by $
 100 ^{n/\varepsilon}$;  a classical problem in convexity is to " effectively enumerate"
 these $P$.

 \section  {Volume Stability In The Riemannian Category.} 

The Burago-Ivanov Volume rigidity theorem (see section 5), like any other sharp inequality between two geometric invariants,
   $$inv_1(X) /inv_2(X)\geq 1,$$
where the equality implies that $X$ is isometric to a particular space $X_0$ (or a member of   a "small explicit"  class of manifolds) raises the following 
 
 \vspace {1mm}

{ \it Stability Problems.}  (a)  Describe spaces  $X$ where 
$inv_1(X)/inv_2 \leq  1+\varepsilon$;

(b)  find a metric in the space  $\cal X$ of all $X$, such that the subspace
${\cal X}_\varepsilon \subset \cal X$ of those $X$ where $inv_1(X)/inv_2 \leq  1+\varepsilon$ is compact.

   \vspace {1mm}

Prior to engaging into general discussion on relevant metrics in the space of metric spaces (we shall come to this in Part 2 of the paper and in \cite {bil} ) it is instructive  
 to look at a  "metric" in the space of Riemannian manifolds, which (almost) adequately reflects the volume rigidity picture.

 \vspace {1mm}

\subsection {\it Directed Lipschitz Metric Normalized by Volume.} To simplify/normalize, let $X$ and $X'$ be  closed connected oriented  Riemannian $n$-manifolds with $vol(X)=vol(X')=V$.   

  Define ${\overset {\longrightarrow}{dist}}_{Lip/vol}(X',X)$ as the infimum of $\varepsilon\geq 0$, such that $X'$ admits an $e^\varepsilon$-Lipschitz
map $L:X'\to X$ of degree $\pm1$, where, observe, $e^\varepsilon=1+\varepsilon +o(\varepsilon)$.

This "metric" is  {\it  non-symmetric},  it may be equal   $+\infty$, but it satisfies the {\it triangle inequality} and it equals zero  if and only if $X$ and $X'$ are isometric. 

\vspace {1mm}

There is nothing wrong with being non-symmetric.   Limits make perfect sense for such "metrics". 

 For example, the  "true", in the category theoretic sense,  Hausdorff metric on subsets $Y\subset X$, call it   $\overset {\longrightarrow}{dist}_{Hau}(Y_1,Y_2)$ --  {\it the minimal $\varepsilon$ such that  $Y_1$ is contained in the $\varepsilon$-neighbourhood of $Y_2$}, is also non-symmetric.

 \vspace {1mm}
 
 Let us  show that suitable bounds on the diameters and curvatures  guaranty an  "almost isometric"
 diffeomorphism  between   $Lip/vol$-close manifolds.

Denote  the minimium of the diameters
of  $X$ and $X'$ by 
$$ diam_{min}=min(diam(X), diam(X'))$$
and  the maximum of
the absolute values of  their sectional curvatures at all tangent $2$-planes  
by$$ |curv|_{max}=max\big (\sup_{\tau_2 \subset T(X)}|sect.curv_{\tau_2}|(X),  \sup_{\tau_2' \subset T(X')}|sect.curv_{\tau_2'}|(X')\big).$$

 \vspace {1mm}

$  \overset {\longrightarrow}{\mathbf  Dist}\leq\varepsilon \Rightarrow {\textbf Diff}$.  {\it   Given numbers $n=1,2,...$, and $c,D\geq 0$, there exists an  $\varepsilon=\varepsilon(n,  cD^2)>0$,   such that the Riemannian $n$-manifolds  $X'$ and $X$ satisfying the inequalities
$$ diam_{min} \leq D, \mbox { }     |curv|_{max}\leq c \mbox { and }\overset {\longrightarrow}{dist}_{Lip/vol}(X',X)\leq \varepsilon$$ 
 are diffeomorphic.
 
 Moreover, the implied diffeomorphism $X'\to X$ 
 is $(1+\delta)$-bi-Lipschitz with $\delta\to 0$ for $\varepsilon \to 0$.}
 
 \vspace {1mm}
 
 To show this,  rescale the metrics in $X$ and $X'$ in order to have $|curv|_{max}=1$ and then approximate the implied $e^\varepsilon$-Lipschitz map $L:X'\to X$ by a $(1+\delta)$-bi-Lipschitz diffeomorphism $\bar L:X'\to X$
 as follows. 
  
  Denote by $\tilde B_x(\rho)\underset {\exp}  \to X$ the $\rho$ ball in the tangent space $T_x(X)$ with the metric induced  from  $X$  by the exponential map and observe that
if   $\rho $ is sufficiently small, $\rho\leq \rho_0( n, D)>0$, and $\varepsilon$ is much smaller than $\rho$, then,  
   for every $x'\in X$, 
   
   \vspace {1mm}
   
   {\it the map $L$ from the ball $B_{x'}(\rho) \subset X'$ to  $B_x(\rho e^\varepsilon) \subset X$, $x=L(x')$,  uniquely lifts to a map   $\tilde L_{x'}:\tilde B_{x'}(\rho) \to \tilde B_x(\rho e^\varepsilon)$.}

Indeed, if otherwise, the map $L$ would send  some short geodesic loops from $X'$   to loops contractible in their immediate vicinity in $X$;  this, in turn, would imply that $deg(L)=0$. 
   
   \vspace {1mm}

Now we turn $L$ into a diffeomorphism as follow. Take  the center of mass $\tilde x \in \tilde B_x(\rho e^\varepsilon) \subset X$   of the Riemannian measure of  $\tilde B_{x'}(\rho)$ mapped to    $\tilde B_x(\rho e^\varepsilon) $  by $\tilde L_{x'} $    and  let $\bar L(x')=\exp_x(\tilde x) \in X$. It is easy to see,   as in \cite {geo-kar}, \cite{kar  },  that  for small $\rho>0$ and $\varepsilon << \rho$,   

  \vspace {1mm}

{\it the so defined   map $\bar L: X'\to X$ is a $(1+\delta)$-Lipschitz diffeomorphism, where $\delta\to 0$ for $\rho \to 0$ and
$\varepsilon/\rho \to 0$.}
   
   \vspace {1mm}

  It is significant that the above does not work if you replace the  bound on the diameters of $X$ and $X'$ by volume bounds, i.e. if you  bound $|curv|_{max}V^{2/n}$ while allowing  $|curv|_{max}\cdot diam_{min}^2 \to \infty$.  

   \vspace {1mm}

{\it (Counter)example}.  Let $X_0$ and $X_0'$  be   complete Riemannian $(n-1)$ manifolds  that are cylindrical at infinity and  have their sectional curvatures bounded by $|sect.curv|\leq 1/2$.
 Let $L_0: X'_0\to X_0$ be a $1$-Lipschitz  map that is isometric  at infinity.

  Thus, both  manifold equal $Y\times [0,\infty)$  on their common cylindrical end for a closed $(n-2)$-manifold $Y$,  where we denote by $t(x)=t(x') \in [0,\infty)$, $x\in X_0$, $x'\in X'_0$, the axial parameter on this end.
  
   (All  this can be arranged, for instance, if $X$ is diffeomorphic to $\mathbb R^{n-1}$ and $X'$ is diffeomorphic to the complement of a ball in an {\it arbitrary} closed $(n-1)$-manifold, where $X_0$ and $X_0'$ can be made isometric to  $S^{n-2}\times [0,\infty)$ at infinity.)

  Let $f$  be a smooth positive function on $X_0$  and  let $f'(x')=(1-\varepsilon)f(L_0(x'))$, $x'\in X'_0$.   (One may take an $f$  which      depends only on $t$ on the cylindrical end of $X_0$, i.e. $f(x)=\phi(t(x))$, and which is constant away from this end.)

  Multiply $X_1$ and $X_1'$  by the unit circle 
  $ S^1$  and 
  modify the product 
  metrics in $X_1= X_0\times S^1$ and 
  $X_1' = X'_0\times S^1$ by multiplying the circles  over 
  $x_0\in X_0$ and $x_0'\in X_0'$, that are $x_0\times S^1$ and $x'_0\times S^1$,
  by the functions $f$ and $f'$ respectively. Denote by   $X_f$ and  $X_{f'}'$
 the resulting complete Riemannain $n$-manifolds, where we organize the matter with a suitable $f$, such that both manifolds  $X_f$ and  $X_{f'}'$ have $|sect.curv| \leq 1$. 

Observe   that the map $L_1 : X'_{f'}=X'_0\times S^1\to X_f = X_0\times S^1$ for $(x',s)\mapsto (L_0(x'),s)$
is $(1+\varepsilon)$-Lipschitz for these metrics.
  
  There obviously exists   an axial  value $t_\varepsilon \in [0, \infty)$, such that
  the integrals of the two functions over the cut-off manifolds $X_0(t_\varepsilon)\subset X_0$ and  $X'_0(t_\varepsilon)\subset X'_0$,  defined by $ t(x), t(x')\leq t_\varepsilon$, are equal. Then the corresponding  compact $n$-manifolds   $X_f(t_\varepsilon)\subset X_f$ and  $X_{f'}'(t_\varepsilon)\subset X_{f'}' $,  with the  boundaries 
  $\partial  X_f(t_\varepsilon)= \partial X'_{f'}(t_\varepsilon)=Y\times S^1 \times t_\varepsilon$,  have equal $n$-volumes.

 If we replace $f\mapsto \epsilon f$ and $f'\mapsto \epsilon f'$ with an arbitrary $\epsilon>0$, we do not change the curvatures of the manifolds; thus we can make the volumes of  $X_f(t_\varepsilon)$ and  $X_{f'}'(t_\varepsilon)$ arbitrarily small, say both equal $\epsilon/2$ (to avoid more  letters in the notation) still keeping $|sect.curv|\leq 1$.

  Finally, we  take the doubles of these manifolds and obtain, for all $n\geq 3$ and arbitrarily small $\varepsilon,\epsilon>0$, lots of
  
  {\it closed  Riemannian $n$-manifolds   
 $X=X(\varepsilon, \epsilon)$ and $X'=X'(\varepsilon, \epsilon)$ both with $|sect.curv|\leq 1$  and with the volumes $=\epsilon$, where $X$ admits a $(1+\varepsilon)$-Lipschitz
 map $L:X'\to X$ of degree $1$ that is not a homotopy equivalence.}
  (The diameters of these manifolds are $\approx t_\varepsilon \approx \varepsilon^{-1}$.)
 
 \vspace {1mm}

On can recapture, however, the implication $  \overset {\longrightarrow}{\mathbf  Dist}\leq\varepsilon \Rightarrow {\textbf Diff}$  for complete Riemannian manifolds with $|sect.curv|\leq 1$ if, for a given $d=1,2,...$,
  the condition $vol(X')=vol(X)$
is replaced  by
$$vol (L^{-1}(B))\leq (d+\varepsilon)vol(B)\mbox { for all unit balls $B\subset X$}. \leqno 
[B'\lesssim_\varepsilon d\cdot B]$$
Then, the  argument used for $  \overset {\longrightarrow}{\mathbf  Dist}\leq\varepsilon \Rightarrow {\textbf Diff}$ also delivers

{\it   an approximation of every  proper $(1+\varepsilon)$-Lipschitz map 
$L: X\to X'$ of  degree $d$ by  a   locally 
diffeomorphic locally  $(1+\delta)$-bi-Lipschitz   map $X'\to X$ for all  sufficiently small $\varepsilon\leq \varepsilon_0(n)>0$ and $\delta\leq \delta(\varepsilon) \underset {\varepsilon \to 0} \to 0$. }
(Such an approximating map  $X'\to X$  necessarily is a $d$-sheeted covering map.)

 \vspace {1mm}

{\it Questions.} Let $X$ and $X'$ be Riemannian $n$-manifolds (or possibly singular Alexandrov spaces for this matter) with their sectional curvatures bounded from below by $-1$ and let $L: X'\to X$ be a   $(1+\varepsilon)$-Lipschitz map of  degree $d$ that satisfy $[B'\lesssim_\varepsilon d\cdot B]$. 

Does then,  for small $\varepsilon >0$, the map $L$ lift  to  a {\it homotopy equivalence} between $X'$ and a $d$-sheeted 
covering of $X$?  (This is easy in the non-collapsed case.)

Is, moreover, $L$ homotopic to a locally homeomorphic map?
 
Is there anything meaningful here with the  lower bound on the sectional curvatures relaxed to a lower bound on the Ricci curvatures of $X$ and/or $X'$?

 \vspace {1mm}
 
 What happens for   $dim(X')-dim(X)=k>0$?
 
 Namely, let  $L:X'\to X$ be an $(1+\varepsilon)$-Lipschitz map such that  $vol(X')\leq ||deg_k(L)||_{vol_k}\cdot vol(X)$, where  $||deg(L)||_{vol_k}$ is the {\it $\mathbb R$-mass} of the homology  class  $deg_k(L) \in H_k(X')$  that is  the  class represented by a  generic pullback  $L^{-1}(x) \subset X'$, $x\in X$. 
 
  Are  there particular  topological/homological constrains on such an $L$ for small $\varepsilon$  in the case where    $X$ and $X' $ satisfy   specified bounds on   their sizes and the local geometries, say an upper bound on the diameters and on the  absolute values of the sectional curvatures?

  \subsection {Mean Curvature Stability.}
  
Below is   another kind of situation where  Lipschitz and volume confront one another.

  Let $X$ and $X'$ be 
 smooth Riemannian $n$-manifolds, $L=L_\varepsilon : X'\to X$  be a $(1+\varepsilon)$-bi-Lipschitz homeomorphism and let 
$\phi :X\to \mathbb R$ be a continuous function.

   Let $ Y\subset X$ be
  a smooth  closed cooriented hypersurface with mean curvature $mn.curv_y(Y)=\phi(y)$ and let positive
   numbers $\rho, \epsilon>0$ be
   given. 

   \vspace {1mm}
  
 If $\varepsilon \leq \varepsilon_0= \varepsilon_0(X, \phi, \rho, \epsilon )>0$, then 
 there exists a closed  hypersurface $Y_{min'}\subset X$ such that the following conditions
 \textbf{[U}$_\rho]$, \textbf {[mn.curv$\pm\epsilon]$} and   \textbf{[diff]} are satisfied.
  
  \vspace {1mm}
 
 \hspace {-5mm}  \textbf{[U}$_\rho]$  \hspace {2mm} $Y_{min'}$ {\it is contained in the
  $\rho$-neighbourhood $U_\rho(Y)\subset X $   of    $Y$
 where it is  homologous to $Y$.}
    
  \vspace {1mm} 
 
 \hspace {-5mm}  \textbf{ [mn.curv$\pm\epsilon]$}  \hspace {2mm}   {\it The hypersurface
  $Y'_{min} =L^{-1}(Y_{min'}) \subset X'$ is $C^2$-smooth and its mean curvature with respect to the Riemannian metric in $X'$ satisfies
$$|mn.curv_{y'}(Y'_{min})-\phi'(y')| \leq \epsilon \mbox { for $\phi'(y')= \phi\circ L(y')$ 
 and   all } y'\in Y'_{min}.$$}
  
  Notice that even if $L$ is smooth, the hypersurface   $Y_{min'} \subset X$ is
  {\it not,  in general,  $C^1$-close
 to} $Y$. Moreover,  even for $dim(Y)=1$, the normal projection $  Y_{min'} \to Y$ 
 is {\it not,} typically, one-to-one. However, 

    \vspace {1mm}

  \hspace {-5mm}     \textbf{[diff]} \hspace {2mm} {\it the manifold $Y'_{min}$ is diffeomorphic to $Y$. In fact the composition of 
$L: Y'_{min} \to  U_\rho(Y)$ with the normal projection  $U_\rho(Y)\to Y$ can be approximated
by a diffeomorphism   $Y'_{min} \to  Y$.}

  \vspace {1mm}

{\it Proof.} We may assume that $Y$ is connected and $X\supset Y$ equals a small normal neighbourhood of $Y$. Thus,   $Y$ divides  $X$ into two halves, call them {\it inside and outsides} of  $Y$, written $in(Y)\subset X $
 and $out(Y) \subset X$, with common boundary $\partial( in(Y))= \partial (out(Y)) =Y$.
 
 \vspace {1mm}
 
 {\it $\mu$-Area and $\mu$-Bubbles.}   Given a measure $\mu$ on $X$,
let 
  $$area_{-\mu}(Y)=_{def} vol_{n-1}(Y)-\mu (in(Y)) $$
 and call a  hypersurface $Y \subset X$ a {\it stable $\mu$-bubble} 
  if it  {\it locally minimizes}  
 the function $Y\mapsto area_{-\mu}(Y)$ among all hypersurfaces in $X$ homologous to $Y$.

If $\mu$ is given by a continous density function
 $\phi(x)$,  $x\in X$, i.e.  $ \mu =\phi \cdot vol_n$ for the Riemannian $n$-volume (measure) 
 $vol_n$, then  these are called {\it $\phi$-bubles}. Clearly, the mean curvature  of a
 $\phi$-bubble $Y \subset X$
satisfies $mn.curv(y)=\phi(y)$. 

In particular, {\it $\varepsilon$-bubbles},  where $\mu$ equals the Riemannian $n$-volume $vol_n$ in $X$ times 
  $\varepsilon \in \mathbb R$, have constant mean curvature $\varepsilon.$

  \vspace {1mm}

 {\it Local Traps.}
 Let $Y\subset X$ be a closed  
smooth cooriented hypersurface and $\phi(x)$ be a $C^1$-smooth function
such that $\phi(y)=mn.curv_y(Y)$ for all $y\in Y$. 

 If {\it the inward normal derivative $\frac {d\phi(y)}{d\nu^{in}}$ on $Y$,
 is sufficiently large}, namely

$$\frac {d\phi(y)}{d\nu^{in}_y}>  curv_y^2( Y)  +Ricci_X(\nu^{in}_y, \nu^{in}_y)\mbox { for all } y \in Y,$$
where $curv^2$ denotes the sum of squares of the principal curvatures of $Y$ and where, observe, $Ricci(\nu^{in}, \nu^{in})=Ricci(\nu^{out}, \nu^{out})$,  then

\vspace {1mm}
{\it $Y$ is a stable $\phi$-bubble; moreover, there is a (small) neighbouthood $U_0\subset X$ of $Y$,  such that every hypersurface $Y'\neq Y \subset U_0$ homologous to $Y$ has strictly 
greater $\phi$-area
than $Y$.}

\vspace {1mm}

 This trivially     follows from  the {\it second variation formula} for $vol_{n-1}(Y)$.

\vspace {1mm}

Now, given a smooth hypesurface $Y\subset X$, let $\phi: X\to \mathbb R$ be equal the mean curvature of $Y$ on $Y$ and have  large  inward normal derivative. Then (almost) obviously, the domain $L_\varepsilon^{-1}(U_0)\subset X'$ contains a locally minimal $\phi'$-bubble, say   $Y'_{min} \subset  L_\varepsilon^{-1}(U_0) $,  homologous to  $L_\varepsilon^{-1}(Y) \subset X'$  for all  $\varepsilon\leq \varepsilon(X, Y, \phi)> 0$

The hypersurface   $Y'_{min}$ can be, a priori, singular. However, if $\varepsilon$ is sufficiently small it is smooth by the following standard argument/

  Let $r>0$ be a small number 
that is, however, 
  much bigger then $\varepsilon$.
Then the volume of
the intersection of  $Y'_{min} $ with the  $r$-ball $B'_y(r) \subset X'$ is
  about $\varepsilon$-close to the volume of the Eucliden
  ball $B_{Eucl}^{n-1}(r)$ for all $y\in Y'_{min}.$  Since $X'$ on the $r$-scale has almost the same
   filling inequalities as $\mathbb R^{n-1}$,
  it follows, by the standard monotonicity argument, that the ratio 
  $$\frac{vol_{n-1}(B'_y(r)\cap Y'_{min} )}{vol_{n-1}(B_{Eucl}^{n-1}(r))}$$
is close to one for all  {\it arbitrarily  small } balls.

Hence, by Almgren-Allard   regularity theory (that is the only  non-elementary ingredient of our argument) the hypersurface $Y'_{min} $
 is $C^2$-smooth with the mean curvature equal $\phi'_\epsilon$ on it.

It remains to show that $Y'_{min} $ is diffeomorphic to $Y$.

You may assume $L$ is smooth, this does not cost you anything for small $\varepsilon$,
but the geometric measure theory tells you nothing,  a priori, about the topology and geometry of
  $Y'_{min}  \subset X'$  and of the corresponding  $Y_{min'}=L(Y'_{min})   \subset X$.   
  
  Yet, 
   since $Y'_{\min}$ is a {\it minimal} bubble,  and  $L$ is $(1+  \varepsilon)$-bi-Lipshitz, the hypersurface is
  $Y_{min'}$
  is  {\it $\varepsilon'$-quasi-minimal} for  $vol_n$ and $\varepsilon'\approx \varepsilon$
    in the sense of \cite {alm}
   on {\it all  sufficiently small scales}. Then   the standard blow-up  rescaling/limit 
   argument yields the following.
  
 {\it Weak Distortion Property.} If the distance function from a point $x\in X$ to  $Y_{min'} $
 has two minima $y_1, y_2 \in Y_{min'} $,
 $$dist (x,y_1)=dist (x,y_2)=\delta=dist (x, Y_{min'} ),$$
 then the angle between the corresponding minimal geodesic
 segments $[x,y_1]$ and $[x,y_2]$ at $x$ must be small, roughly of order $\varepsilon$, 
 for all not very large 
$ \delta$.
  
 It follows that the distance function $x\to dist(x,Y_{min'})$, say outward $Y_{min'}$,  can be
  "bi-Lipschitz" approximated in  the $\delta$-neighbourhood $ U_\delta(Y_{min})$ for, say 
  $\delta=100^n\rho$ for a sufficiently small $\rho$ and $\epsilon$.
   by a smooth function $d(x)$ {\it  without critical points } which vanishes on $Y_{min'}$. 
   Similarly one sees
   that the normal projection of the $\delta/2$-level
   of $d(x)$ to $Y$ is a diffeomorphism.  QED. 
   
    \vspace {2mm}
    
{\it Concluding Remarks.}     The $(1+\varepsilon)$-bi-Lipschitz assumption, as well
as the $e^\varepsilon$-Lipschitz constrain on maps in the previous section, are unduly restrictive.
These will be relaxed in the Hilbert volume framework, in the spirit of Stephan Wenger's metric in the space of manifolds, 
\cite {wen}  \cite {wen1}, in Part 2 of our paper.

  \end{document}